\input amstex

\def\next{AMS-P}

\def\next{AMSPPT}
\ifx\styname\next \else\input amsppt.sty \relax\fi



\brokenpenalty=10000
\clubpenalty=10000
\widowpenalty=10000

\catcode`\@=11

\def\keyboarder#1{}%

\def\pagewidth#1{\hsize#1
  \captionwidth@24pc}

\pagewidth{30pc}

\parindent=18\p@
\normalparindent\parindent

\parskip=\z@

\def\foliofont@{\sevenrm}
\def\headlinefont@{\sevenpoint}

\def\leftheadline{\rlap{\foliofont@\folio}\hfill \iftrue\topmark\fi \hfill}
\def\rightheadline{\hfill \expandafter
  \hfill \llap{\foliofont@\folio}}

\def\specialheadfont@{\elevenpoint\smc}
\def\headfont@{\bf}
\def\subheadfont@{\bf}
\def\refsheadfont@{\bf}
\def\abstractfont@{\smc}
\def\proclaimheadfont@{\smc}
\def\xcaheadfont@{\smc}
\def\captionfont@{\smc}
\def\citefont@{\bf}
\def\refsfont@{\eightpoint}

\font\sixsy=cmsy6


\font@\twelverm=cmr10 scaled \magstep1
\font@\twelvebf=cmbx10 scaled \magstep1
\font@\twelveit=cmti10 scaled \magstep1
\font@\twelvesl=cmsl10 scaled \magstep1
\font@\twelvesmc=cmcsc10 scaled \magstep1
\font@\twelvett=cmtt10 scaled \magstep1
\font@\twelvei=cmmi10 scaled \magstep1
\font@\twelvesy=cmsy10 scaled \magstep1
\font@\twelveex=cmex10 scaled \magstep1
\font@\twelvemsb=msbm10 scaled \magstep1
\font@\twelveeufm=eufm10 scaled \magstep1

\newtoks\twelvepoint@
\def\twelvepoint{\normalbaselineskip14\p@
  \abovedisplayskip12\p@ plus3\p@ minus9\p@
  \belowdisplayskip\abovedisplayskip
  \abovedisplayshortskip\z@ plus3\p@
  \belowdisplayshortskip7\p@ plus3\p@ minus4\p@
  \textonlyfont@\rm\twelverm \textonlyfont@\it\twelveit
  \textonlyfont@\sl\twelvesl \textonlyfont@\bf\twelvebf
  \textonlyfont@\smc\twelvesmc \textonlyfont@\tt\twelvett
  \ifsyntax@ \def\big##1{{\hbox{$\left##1\right.$}}}%
    \let\Big\big \let\bigg\big \let\Bigg\big
  \else
    \textfont\z@\twelverm  \scriptfont\z@\eightrm
       \scriptscriptfont\z@\sixrm
    \textfont\@ne\twelvei  \scriptfont\@ne\eighti
       \scriptscriptfont\@ne\sixi
    \textfont\tw@\twelvesy \scriptfont\tw@\eightsy
       \scriptscriptfont\tw@\sixsy
    \textfont\thr@@\twelveex \scriptfont\thr@@\eightex
        \scriptscriptfont\thr@@\eightex
    \textfont\itfam\twelveit \scriptfont\itfam\eightit
        \scriptscriptfont\itfam\eightit
    \textfont\bffam\twelvebf \scriptfont\bffam\eightbf
        \scriptscriptfont\bffam\sixbf
    \textfont\msbfam\twelvemsb \scriptfont\msbfam\eightmsb
        \scriptscriptfont\msbfam\sixmsb
    \textfont\eufmfam\twelveeufm \scriptfont\eufmfam\eighteufm
        \scriptscriptfont\eufmfam\sixeufm
    \setbox\strutbox\hbox{\vrule height8.5\p@ depth3.5\p@ width\z@}%
    \setbox\strutbox@\hbox{\lower.5\normallineskiplimit\vbox{%
        \kern-\normallineskiplimit\copy\strutbox}}%
    \setbox\z@\vbox{\hbox{$($}\kern\z@}\bigsize@1.2\ht\z@
  \fi
  \normalbaselines\rm\dotsspace@1.5mu\ex@.2326ex\jot3\ex@
  \the\twelvepoint@}

\font@\elevenrm=cmr10 scaled \magstephalf
\font@\elevenbf=cmbx10 scaled \magstephalf
\font@\elevenit=cmti10 scaled \magstephalf
\font@\elevensl=cmsl10 scaled \magstephalf
\font@\elevensmc=cmcsc10 scaled \magstephalf
\font@\eleventt=cmtt10 scaled \magstephalf
\font@\eleveni=cmmi10 scaled \magstephalf
\font@\elevensy=cmsy10 scaled \magstephalf
\font@\elevenex=cmex10 scaled \magstephalf
\font@\elevenmsb=msbm10 scaled \magstephalf
\font@\eleveneufm=eufm10 scaled \magstephalf

\newtoks\elevenpoint@
\def\elevenpoint{\normalbaselineskip13\p@
  \abovedisplayskip12\p@ plus3\p@ minus9\p@
  \belowdisplayskip\abovedisplayskip
  \abovedisplayshortskip\z@ plus3\p@
  \belowdisplayshortskip7\p@ plus3\p@ minus4\p@
  \textonlyfont@\rm\elevenrm \textonlyfont@\it\elevenit
  \textonlyfont@\sl\elevensl \textonlyfont@\bf\elevenbf
  \textonlyfont@\smc\elevensmc \textonlyfont@\tt\eleventt
  \ifsyntax@ \def\big##1{{\hbox{$\left##1\right.$}}}%
    \let\Big\big \let\bigg\big \let\Bigg\big
  \else
    \textfont\z@\elevenrm  \scriptfont\z@\eightrm
       \scriptscriptfont\z@\sixrm
    \textfont\@ne\eleveni  \scriptfont\@ne\eighti
       \scriptscriptfont\@ne\sixi
    \textfont\tw@\elevensy \scriptfont\tw@\eightsy
       \scriptscriptfont\tw@\sixsy
    \textfont\thr@@\elevenex \scriptfont\thr@@\eightex
        \scriptscriptfont\thr@@\eightex
    \textfont\itfam\elevenit \scriptfont\itfam\eightit
        \scriptscriptfont\itfam\eightit
    \textfont\bffam\elevenbf \scriptfont\bffam\eightbf
        \scriptscriptfont\bffam\sixbf
    \textfont\msbfam\elevenmsb \scriptfont\msbfam\eightmsb
        \scriptscriptfont\msbfam\sixmsb
    \textfont\eufmfam\eleveneufm \scriptfont\eufmfam\eighteufm
        \scriptscriptfont\eufmfam\sixeufm
    \setbox\strutbox\hbox{\vrule height8.5\p@ depth3.5\p@ width\z@}%
    \setbox\strutbox@\hbox{\lower.5\normallineskiplimit\vbox{%
        \kern-\normallineskiplimit\copy\strutbox}}%
    \setbox\z@\vbox{\hbox{$($}\kern\z@}\bigsize@1.2\ht\z@
  \fi
  \normalbaselines\rm\dotsspace@1.5mu\ex@.2326ex\jot3\ex@
  \the\elevenpoint@}

\addto\tenpoint{\normalbaselineskip12\p@
 \abovedisplayskip6\p@ plus6\p@ minus0\p@
 \belowdisplayskip6\p@ plus6\p@ minus0\p@
 \abovedisplayshortskip0\p@ plus3\p@ minus0\p@
 \belowdisplayshortskip2\p@ plus3\p@ minus0\p@
 \ifsyntax@
 \else
  \setbox\strutbox\hbox{\vrule height9\p@ depth4\p@ width\z@}%
  \setbox\strutbox@\hbox{\vrule height8\p@ depth3\p@ width\z@}%
 \fi
 \normalbaselines\rm}

\newtoks\sevenpoint@
\def\sevenpoint{\normalbaselineskip9\p@
 \textonlyfont@\rm\sevenrm \textonlyfont@\it\sevenit
 \textonlyfont@\sl\sevensl \textonlyfont@\bf\sevenbf
 \textonlyfont@\smc\sevensmc \textonlyfont@\tt\seventt
  \textfont\z@\sevenrm \scriptfont\z@\sixrm
       \scriptscriptfont\z@\fiverm
  \textfont\@ne\seveni \scriptfont\@ne\sixi
       \scriptscriptfont\@ne\fivei
  \textfont\tw@\sevensy \scriptfont\tw@\sixsy
       \scriptscriptfont\tw@\fivesy
  \textfont\thr@@\sevenex \scriptfont\thr@@\sevenex
   \scriptscriptfont\thr@@\sevenex
  \textfont\itfam\sevenit \scriptfont\itfam\sevenit
   \scriptscriptfont\itfam\sevenit
  \textfont\bffam\sevenbf \scriptfont\bffam\sixbf
   \scriptscriptfont\bffam\fivebf
  \textfont\msbfam\sevenmsb \scriptfont\msbfam\sixmsb
   \scriptscriptfont\msbfam\fivemsb
  \textfont\eufmfam\seveneufm \scriptfont\eufmfam\sixeufm
   \scriptscriptfont\eufmfam\fiveeufm
 \setbox\strutbox\hbox{\vrule height7\p@ depth3\p@ width\z@}%
 \setbox\strutbox@\hbox{\raise.5\normallineskiplimit\vbox{%
   \kern-\normallineskiplimit\copy\strutbox}}%
 \setbox\z@\vbox{\hbox{$($}\kern\z@}\bigsize@1.2\ht\z@
 \normalbaselines\sevenrm\dotsspace@1.5mu\ex@.2326ex\jot3\ex@
 \the\sevenpoint@}

\def\loadeusm{\loadmathfont{eusm}%
  \addto\sevenpoint{\textfont\eusmfam\seveneusm \scriptfont\eusmfam\sixeusm
    \scriptscriptfont\eusmfam\fiveeusm}%
  \font@\eleveneusm=eusm10 scaled \magstephalf
  \addto\elevenpoint{\textfont\eusmfam\eleveneusm \scriptfont\eusmfam\eighteusm
    \scriptscriptfont\eusmfam\sixeusm}%
  \font@\twelveeusm=eusm10 scaled \magstep1
  \addto\twelvepoint{\textfont\eusmfam\twelveeusm \scriptfont\eusmfam\eighteusm
    \scriptscriptfont\eusmfam\sixeusm}%
  }

\def\xci@{1991}
\def\mm@{2000}
\def\subjyear@{1991}
\def\subjclassyear#1{%
  \def\next{#1}%
  \ifx\next\mm@ \def\subjyear@{#1}%
  \else\ifx\next\xci@ \def\subjyear@{#1 }%
  \else \err@{Unknown edition (#1) of Mathematics
      Subject Classification; using 1991'.}
  \fi\fi}

\newskip\abovespecheadskip   \abovespecheadskip20\p@ plus8\p@ minus2\p@
\newdimen\belowspecheadskip  \belowspecheadskip6\p@
\outer\def\specialhead{%
  \add@missing\endroster \add@missing\enddefinition
  \add@missing\enddemo \add@missing\endexample
  \add@missing\endproclaim
  \penaltyandskip@{-200}\abovespecheadskip
  \begingroup\interlinepenalty\@M\rightskip\z@ plus\hsize
  \let\\\linebreak
  \specialheadfont@\raggedcenter@\noindent}

\let\varindent@\indent

\let\subsubhead\relax
\outer\def\subsubhead{%
  \add@missing\endroster \add@missing\enddefinition
  \add@missing\enddemo
  \add@missing\endexample \add@missing\endproclaim
  \let\savedef@\subsubhead \let\subsubhead\relax
  \def\subsubhead##1\endsubsubhead{\restoredef@\subsubhead
    {\def\usualspace{\/{\subsubheadfont@\enspace}}%
     \subsubheadfont@##1\unskip\frills@{.\enspace}}\ignorespaces}%
  \nofrillscheck\subsubhead}


\newskip\abstractindent 	\abstractindent=3pc
\long\def\block #1\endblock{\vskip 6pt
	{\leftskip=\abstractindent \rightskip=\abstractindent
	\noindent #1\endgraf}\vskip 6pt}

\long\def\ext #1\endext{\removelastskip\block #1\endblock}

\outer\def\xca{\let\savedef@\xca \let\xca\relax
  \add@missing\endproclaim \add@missing\endroster
  \add@missing\endxca \envir@stack\endxca
  \def\xca##1{\restoredef@\xca
    \penaltyandskip@{-100}\medskipamount
    \bgroup{\def\usualspace{{\xcaheadfont@\enspace}}%
      \varindent@\xcaheadfont@\ignorespaces##1\unskip
      \frills@{.\xcaheadfont@\enspace}}%
      \ignorespaces}%
  \nofrillscheck\xca}
\def\endxca{\egroup\revert@envir\endxca
  \par\medskip}

\def\remarkheadfont@{\smc}
\def\remark{\let\savedef@\remark \let\remark\relax
  \add@missing\endroster \add@missing\endproclaim
  \envir@stack\endremark
  \def\remark##1{\restoredef@\remark
    \penaltyandskip@{-100}\medskipamount
    {\def\usualspace{{\remarkheadfont@\enspace}}%
     \varindent@\remarkheadfont@\ignorespaces##1\unskip
     \frills@{.\enspace}}\rm
    \ignorespaces}\nofrillscheck\remark}
\def\endremark{\par\revert@envir\endremark\medskip}

\def\qed{\ifhmode\unskip\nobreak\fi\hfill
  \ifmmode\square\else$\m@th\square$\fi}
\def\qedhere{\global\@qedheretrue
  \ifmmode\qquad\square\else\unskip\nobreak\hfill$\m@th\square$\fi}


\newdimen\rosteritemsep
\rosteritemsep=.5pc

\newdimen\rosteritemwd@
\newdimen\rosteritemitemwd
\newdimen\rosteritemitemitemwd

\newbox\setwdbox
\setbox\setwdbox\hbox{(0)}\rosteritemwd=\wd\setwdbox
\rosteritemwd@=\rosteritemwd
\setbox\setwdbox\hbox{(0)\hskip.5pc(c)}\rosteritemitemwd=\wd\setwdbox
\setbox\setwdbox\hbox{(0)\hskip.5pc(c)\hskip.5pc(iii)}%
  \rosteritemitemitemwd=\wd\setwdbox

\def\roster{%
  \envir@stack\endroster
  \edef\leftskip@{\leftskip\the\leftskip}%
  \relaxnext@
  \rostercount@\z@
  \def\item{\FN@\rosteritem@}%
  \def\itemitem{\FN@\rosteritemitem@}%
  \def\itemitemitem{\FN@\rosteritemitemitem@}%
  \DN@{\ifx\next\runinitem\let\next@\nextii@
    \else\let\next@\nextiii@
    \fi\next@}%
  \DNii@\runinitem
    {\unskip
     \DN@{\ifx\next[\let\next@\nextii@
       \else\ifx\next"\let\next@\nextiii@\else\let\next@\nextiv@\fi
       \fi\next@}%
     \DNii@[####1]{\rostercount@####1\relax
       \therosteritem{\number\rostercount@}~\ignorespaces}%
     \def\nextiii@"####1"{{\rm####1}~\ignorespaces}%
     \def\nextiv@{\therosteritem1\rostercount@\@ne~}%
     \par@\firstitem@false
     \FN@\next@}
  \def\nextiii@{\par\par@
    \penalty\@m\vskip-\parskip
    \firstitem@true}%
  \FN@\next@}

\def\rosteritem@{\iffirstitem@\firstitem@false
  \else\par\vskip-\parskip
  \fi
  \leftskip\rosteritemwd \advance\leftskip\normalparindent
  \advance\leftskip.5pc \noindent
  \DNii@[##1]{\rostercount@##1\relax\itembox@}%
  \def\nextiii@"##1"{\def\therosteritem@{\rm##1}\itembox@}%
  \def\nextiv@{\advance\rostercount@\@ne\itembox@}%
  \def\therosteritem@{\therosteritem{\number\rostercount@}}%
  \ifx\next[\let\next@\nextii@
  \else\ifx\next"\let\next@\nextiii@\else\let\next@\nextiv@\fi
  \fi\next@}

\def\itembox@{\llap{\hbox to\rosteritemwd{\hss
  \kern\z@ 
  \therosteritem@}\hskip.5pc}\ignorespaces}

\def\rosteritemitem@{\iffirstitem@\firstitem@false
  \else\par\vskip-\parskip
  \fi
  \leftskip\rosteritemitemwd \advance\leftskip\normalparindent
  \advance\leftskip.5pc \noindent
  \DNii@[##1]{\rostercount@##1\relax\itemitembox@}%
  \def\nextiii@"##1"{\def\therosteritemitem@{\rm##1}\itemitembox@}%
  \def\nextiv@{\advance\rostercount@\@ne\itemitembox@}%
  \def\therosteritemitem@{\therosteritemitem{\number\rostercount@}}%
  \ifx\next[\let\next@\nextii@
  \else\ifx\next"\let\next@\nextiii@\else\let\next@\nextiv@\fi
  \fi\next@}

\def\itemitembox@{\llap{\hbox to\rosteritemitemwd{\hss
  \kern\z@ 
  \therosteritemitem@}\hskip.5pc}\ignorespaces}

\def\therosteritemitem#1{\rom{(\ignorespaces#1\unskip)}}

\def\rosteritemitemitem@{\iffirstitem@\firstitem@false
  \else\par\vskip-\parskip
  \fi
  \leftskip\rosteritemitemitemwd \advance\leftskip\normalparindent
  \advance\leftskip.5pc \noindent
  \DNii@[##1]{\rostercount@##1\relax\itemitemitembox@}%
  \def\nextiii@"##1"{\def\therosteritemitemitem@{\rm##1}\itemitemitembox@}%
  \def\nextiv@{\advance\rostercount@\@ne\itemitemitembox@}%
  \def\therosteritemitemitem@{\therosteritemitemitem{\number\rostercount@}}%
  \ifx\next[\let\next@\nextii@
  \else\ifx\next"\let\next@\nextiii@\else\let\next@\nextiv@\fi
  \fi\next@}

\def\itemitemitembox@{\llap{\hbox to\rosteritemitemitemwd{\hss
  \kern\z@ 
  \therosteritemitemitem@}\hskip.5pc}\ignorespaces}

\def\therosteritemitemitem#1{\rom{(\ignorespaces#1\unskip)}}

\def\endroster{\relaxnext@
  \revert@envir\endroster 
  \par\leftskip@ 
  \global\rosteritemwd\rosteritemwd@ 
  \penalty-50 
  \DN@{\ifx\next\Runinitem\let\next@\relax
    \else\nextRunin@false\let\item\plainitem@ 
      \ifx\next\par 
        \DN@\par{\everypar\expandafter{\the\everypartoks@}}%
      \else 
        \DN@{\noindent\everypar\expandafter{\the\everypartoks@}}%
      \fi
    \fi\next@}%
  \FN@\next@}


\def\address#1\endaddress{\global\advance\addresscount@\@ne
  \expandafter\gdef\csname address\number\addresscount@\endcsname
  {\vskip12\p@ minus6\p@\indent\addressfont@\smc\ignorespaces#1\par}}


\def\curraddr{\let\savedef@\curraddr
  \def\curraddr##1\endcurraddr{\let\curraddr\savedef@
    \if\notempty{##1}%
      \toks@\expandafter\expandafter\expandafter{%
        \csname address\number\addresscount@\endcsname}%
      \toks@@{##1}%
      \expandafter\xdef\csname address\number\addresscount@\endcsname
        {\the\toks@\endgraf\noexpand\nobreak
          \indent\noexpand\addressfont@{\noexpand\rm
          \frills@{{\noexpand\it Current address\noexpand\/}:\space}%
          \def\noexpand\usualspace{\space}\the\toks@@\unskip}}%
    \fi}%
  \nofrillscheck\curraddr}


\def\email{\let\savedef@\email
  \def\email##1\endemail{\let\email\savedef@
    \if\notempty{##1}%
      \toks@{\def\usualspace{{\it\enspace}}\endgraf\indent\addressfont@}%
      \toks@@{{\tt##1}\par}%
      \expandafter\xdef\csname email\number\addresscount@\endcsname
      {\the\toks@\frills@{{\noexpand\it E-mail address\noexpand\/}:%
        \noexpand\enspace}\the\toks@@}%
    \fi}%
  \nofrillscheck\email}

\def\bysame{\by\hbox to2pc{\hrulefill}\thinspace\kern\z@}

\def\refstyle#1{\uppercase{%
  \gdef\refstyle@{#1}%
  \if#1A\relax \def\keyformat##1{[##1]\enspace\hfil}%
  \else\if#1B\relax
    \refindentwd2pc
    \def\keyformat##1{\aftergroup\kern
              \aftergroup-\aftergroup\refindentwd}%
  \else\if#1C\relax
    \def\keyformat##1{\hfil##1.\enspace}%
  \fi\fi\fi}
}

\refstyle{A}

\catcode`\@=11

\def\pretitle{\null\vskip74pt}

\def\addressfont@{\eightpoint}


%
%

\define\issueinfo#1#2#3#4{%
  \def\issuevol@{#1}\def\issueno@{#2}%
  \def\issuemonth@{#3}\def\issueyear@{#4}}

\define\originfo#1#2#3#4{\def\origvol@{#1}\def\origno@{#2}%
  \def\origmonth@{#3}\def\origyear@{#4}}

\define\copyrightinfo#1#2{\def\cryear@{#1}\def\crholder@{#2}}

\define\pagespan#1#2{\pageno=#1\def\start@page{#1}\def\end@page{#2}}

\issueinfo{00}{0}{}{1997}
\originfo{00}{0}{}{1997}
\copyrightinfo{\issueyear@}{American Mathematical Society}
\pagespan{000}{000}
\pageno=1 


\def\nojourlogo{\let\jourlogo\empty@}

\def\journame{AMS Proceedings Style}
\def\volinfo{Volume {\sixbf\issuevol@}, \issueyear@}
\let\jourlogoextra@\empty@
\let\jourlogoright@\empty@

\def\jourlogofont@{\sixrm\baselineskip7\p@\relax}
\def\jourlogo{%
  \vbox to\headlineheight{%
    \parshape\z@ \leftskip\z@ \rightskip\z@
    \parfillskip\z@ plus1fil\relax
    \jourlogofont@ \frenchspacing
    \line{\vtop{\parindent\z@ \hsize=.5\hsize
      \journame\newline\volinfo\jourlogoextra@\endgraf\vss}%
      \hfil
      \jourlogoright@
    }%
    \vss}%
}

\def\issn#1{\gdef\theissn{#1}}
\issn{0000-0000}


\def\copyrightline@{%
  \rightline{\sixrm \textfont2=\sixsy \copyright\cryear@\ \crholder@}}

\def\logo@{\copyrightline@}

\def\titlefont{%
 \ifsyntax@\else \twelvepoint\bf \fi }

\def\authorfont{%
  \ifsyntax@
  \else \elevenpoint
  \fi}

\def\title{\let\savedef@\title
  \def\title##1\endtitle{\let\title\savedef@
    \global\setbox\titlebox@\vtop{\titlefont\bf
      \raggedcenter@\frills@{##1}\endgraf}%
    \ifmonograph@
      \edef\next{\the\leftheadtoks}\ifx\next\empty \leftheadtext{##1}\fi
    \fi
    \edef\next{\the\rightheadtoks}\ifx\next\empty \rightheadtext{##1}\fi
  }%
  \nofrillscheck\title}

\def\author#1\endauthor{\global\setbox\authorbox@
  \vbox{\authorfont\raggedcenter@
    {\ignorespaces#1\endgraf}}\relaxnext@
  \edef\next{\the\leftheadtoks}%
  \ifx\next\empty\expandafter\uppercase{\leftheadtext{#1}}\fi}

\def\abstract{\let\savedef@\abstract
  \def\abstract{\let\abstract\savedef@
    \setbox\abstractbox@\vbox\bgroup\indenti=3pc\noindent$$\vbox\bgroup
      \def\envir@end{\endabstract}\advance\hsize-2\indenti
      \def\usualspace{\enspace}\eightpoint \noindent
      \frills@{{\abstractfont@ Abstract.\enspace}}}%
  \nofrillscheck\abstract}

\def\dedicatory #1\enddedicatory{\def\preabstract{{\vskip 20\p@
  \eightpoint\it \raggedcenter@#1\endgraf}}}

\outer\def\endtopmatter{\add@missing\endabstract
  \edef\next{\the\leftheadtoks}%
  \ifx\next\empty@
    \expandafter\leftheadtext\expandafter{\the\rightheadtoks}%
  \fi
  \ifx\thesubjclass@\empty@\else \makefootnote@{}{\thesubjclass@}\fi
  \ifx\thekeywords@\empty@\else \makefootnote@{}{\thekeywords@}\fi
  \ifx\thethanks@\empty@\else \makefootnote@{}{\thethanks@}\fi
  \inslogo@
  \pretitle
  \box\titlebox@
  \topskip10pt
  \preauthor
  \ifvoid\authorbox@\else \vskip16\p@ plus6\p@ minus0\p@\unvbox\authorbox@\fi
  \predate
  \ifx\thedate@\empty\else \vskip6\p@ plus2\p@ minus0\p@
    \line{\hfil\thedate@\hfil}\fi
  \setabstract@
  \nobreak
  \ifvoid\tocbox@\else\vskip1.5pc plus.5pc \unvbox\tocbox@\fi
  \prepaper
  \vskip36\p@\tenpoint
}

\def\setabstract@{%
  \preabstract
  \ifvoid\abstractbox@\else \vskip20\p@ \unvbox\abstractbox@ \fi
}



\begingroup
\let\specialhead\relax
\let\head\relax
\let\subhead\relax
\let\subsubhead\relax
\let\title\relax
\let\chapter\relax

\gdef\setwidest@#1#2{%
   \ifx#1\head\setbox\tocheadbox@\hbox{#2.\enspace}%
   \else\ifx#1\subhead\setbox\tocsubheadbox@\hbox{#2.\enspace}%
   \else\ifx#1\subsubhead\setbox\tocsubheadbox@\hbox{#2.\enspace}%
   \else\ifx#1\key
       \if C\refstyle@ \else\refstyle A\fi
       \setboxz@h{\refsfont@\keyformat{#2}}%
       \refindentwd\wd\z@
   \else\ifx#1\no\refstyle C%
       \setboxz@h{\refsfont@\keyformat{#2}}%
       \refindentwd\wd\z@
   \else\ifx#1\page\setbox\z@\hbox{\quad\bf#2}%
       \pagenumwd\wd\z@
   \else\ifx#1\item
       \setboxz@h{#2}\rosteritemwd=\wd\z@
   \else\ifx#1\itemitem
       \setboxz@h{#2}\rosteritemitemwd=\wd\z@
	\advance\rosteritemitemwd by .5pc
	\advance\rosteritemitemwd by \rosteritemwd
   \else\ifx#1\itemitemitem
       \setboxz@h{#2}\rosteritemitemitemwd=\wd\z@
	\advance\rosteritemitemitemwd by .5pc
	\advance\rosteritemitemitemwd by \rosteritemitemwd
   \else\message{\string\widestnumber\space not defined for this
      option (\string#1)}%
\fi\fi\fi\fi\fi\fi\fi\fi\fi}

\refstyle{A}
\widestnumber\key{M}

\endgroup

\catcode`\@=13

\catcode`\@=11

\def\journame{Contemporary Mathematics}

\issn{0271-4132}

\catcode`\@=13


\NoBlackBoxes
\TagsAsMath

\loadbold

\input xy \xyoption{matrix} \xyoption{arrow}
          \xyoption{curve}  \xyoption{frame}
\def\edge{\ar@{-}}
\def\dttdar{\ar@{.>}}
\def\drbl{\save+<0ex,-2ex> \drop{\bullet} \restore}

\def\dashedge{\ar@{--}}

\def\dropdown#1{\save+<0ex,-4ex> \drop{#1} \restore}
\def\dropup#1{\save+<0ex,4ex> \drop{#1} \restore}

\def\AA{{\Bbb A}}
\def\BB{{\Bbb B}}

\def\PP{{\Bbb P}}
\def\SS{{\Bbb S}}
\def\ZZ{{\Bbb Z}}

\def\RR{{\Bbb R}}

\def\seq{\mathrel{\widehat{=}}}
\def\la{{\Lambda}}
\def\lamod{A\text{-}\roman{mod}}

\def\Amod{A\text{-}\roman{mod}}

\def \len{\operatorname{length}} 

\def\aut{\operatorname{Aut}}

\def\GL{{\operatorname{GL}}}

\def\autlap{\operatorname{Aut}_A(P)}

\def\ext{\operatorname{Ext}}

\def\D{{\Cal D}}

\def\J{{\Cal J}}

\def\C{{\frak C}}
\def\frakD{{\frak D}}

\def\U{{\frak U}}

\def\S{{\sigma}}
\def\vareps{\varepsilon}

\def\bc{{\bold{c}}}
\def\bd{{\bold{d}}}

\def\bP{\bold P}
\def\bPhat{\widehat{\bP}}

\def\rhohat{\hat{\rho}}

\def\Phat{\widehat{P}}

\def\Rhat{\widehat{R}}

\def\hatY{\hat{Y}}

\def\Rep{\operatorname{Rep}}
\def\replad{\Rep(A,d)}
\def\replabd{\Rep(A,\bd)}

\def\repSS{\Rep(\SS)}  
\def\reptbd{\Rep^T_{\bd}}

\def\Schu{\operatorname{\Schu}}

\def\Schubert{\operatorname{\text{\smc{schubert}}}}
\def\SchubertS{\Schubert(\S)}

\def\grasstd{\operatorname{\frak{Grass}}^T_d}
\def\grasstbd{\operatorname{\frak{Grass}}^T_{\bold d}}
\def\GRASS{\operatorname{\text{\smc{grass}}}}

\def\biggrasstbd{\GRASS^T_{\bold d}}

\def\grassbd{\GRASS(A,\bd)}

\def\grassS{\operatorname{\frak{Grass}}(\S)}

\def\grassSS{\operatorname{\frak{Grass}}(\SS)}
\def\biggrassSS{\GRASS(\SS)}

\def\biggrassS{\GRASS(\S)}

\def\grass{\operatorname{\frak{Grass}}}
\def\biggrass{\GRASS}

\def\Hom{\operatorname{Hom}}

\def\Gr{\operatorname{Gr}}

\def\autlabp{\aut_A(\bold{P})}

\def\underbardim{\operatorname{\underline{dim}}}

\def\GLbd{\GL(\bd)}
\def\cbb{c_{\,b',b}}

\def\vsubseteq{\hbox{$\bigcup$\kern0.1em\raise0.05ex\hbox{$\tsize|$}}}
\def\smvsubseteq{\hbox{$\ssize\bigcup$\kern0.03em\raise0.05ex\hbox{$\ssize|$}}}

\def\BHT{{\bf 1}}
\def\codes{{\bf 2}}
\def\GeomIV{{\bf 3}}
\def\CBS{{\bf 4}}
\def\djong{{\bf 5}}
\def\GP{{\bf 6}}
\def\EGAfourtwo{{\bf 7}}
\def\Har{{\bf 8}}
\def\Hess{{\bf 9}}
\def\degen{{\bf 10}}
\def\hierarchies{{\bf 11}}
\def\Kone{{\bf 12}}
\def\Ktwo{{\bf 13}}
\def\Kraft{{\bf 14}}
\def\Mats{{\bf 15}}
\def\MS{{\bf 16}}
\def\Rie{{\bf 17}}
\def\RRS{{\bf 18}}
\def\Scho{{\bf 19}}
\def\Schr{{\bf 20}}
\def\zitat#1{[#1]}

\topmatter

\title Irreducible components of module varieties:  projective equations and rationality
\endtitle

\rightheadtext{Irreducible components of module varieties}

\author  B. Huisgen-Zimmermann and K.R. Goodearl
\endauthor

\address Department of Mathematics, University of California, Santa Barbara, CA 93106
\endaddress

\email birge\@math.ucsb.edu; goodearl\@math.ucsb.edu
\endemail

\abstract     
We expand the existing arsenal of methods for exploring the irreducible components of the varieties $\operatorname{Rep}(A,\bold d)$ which parametrize the representations with dimension vector $\bold d$ of a finite dimensional algebra $A$.  To do so, we move back and forth between $\operatorname{Rep}(A,\bold d)$ and a projective variety, $\operatorname{\text{\smc{grass}}}(A,\bold d)$, parametrizing the same set of isomorphism classes of modules.  In particular, we show the irreducible components to be accessible in a highly compressed format within the projective setting.  Our results include necessary and sufficient conditions for unirationality, smoothness, and normality, followed by applications.  Moreover, they provide equational access to the irreducible components of $\operatorname{\text{\smc{grass}}}(A,\bold d)$ and techniques for deriving qualitative information regarding both the affine and projective scenarios. 
\endabstract

\thanks Both authors were partially supported by grants from the National Science Foundation. \endthanks  

\subjclassyear{2000} \subjclass 16G10, 16G20
\endsubjclass

\endtopmatter

\document

\head 1. Introduction and conventions \endhead

Given a nonhereditary basic finite dimensional algebra $A$ over an 
algebraically closed field $K$, it is typically very difficult  to 
determine/classify/analyze the irreducible components of the affine 
varieties $\replabd$ which parametrize the left $A$-modules with fixed 
dimension vector $\bd$.   Documented interest in this problem dates back 
to 1980, when H.-P. Kraft drew attention to it in his Puebla lectures 
\zitat{\Kraft}.  The first challenging class of algebras for which this 
problem was completely resolved, by Schr\"oer  in \zitat{\Schr}, is a 
subclass of the algebras now dubbed biserial algebras (this subclass was 
rendered prominent by Gelfand and Ponomarev \zitat{\GP}, who used it  as a 
vehicle to understand the representation theory of the Lorentz group).  A 
first installment of  methods designed to tackle the general case was 
assembled in \zitat{\CBS}.  Combining these with the work in \zitat{\Schr}, 
one moreover obtains a classification of the irreducible components of the 
varieties $\replabd$ for the closely related biserial algebra 
$K[X,Y]/(X^2, Y^2)$,  singled out by Carlson due to properties of its 
$\GL(\bd)$-orbit closures; see \zitat{\Rie, Section 3.1} and Section 4 
below. Recently, the components arising in the Carlson example were 
explored from scratch  in \zitat{\RRS};  more detail is given below.  To 
return to the general case:  In \zitat{\BHT}, Babson, Thomas and the first 
author employed alternative projective parametrizing varieties  --  we label 
them  $\grassbd$  --  to build general methodology for determining and 
studying the irreducible components of $\replabd$ (which naturally 
correspond to the components of $\grassbd$), and for analyzing the generic 
behavior of the modules represented by them.   The latter line was 
initiated by Kac for hereditary algebras in his seminal papers 
\zitat{\Kone}, \zitat{\Ktwo}, and continued by Schofield in \zitat{\Scho}; in 
the hereditary  case, the $\replabd$ are full affine spaces, and hence 
irreducible.

An equational study of the components of the varieties $\replabd$ for the 
case of the Carlson algebra was undertaken by Riedtmann, Rutscho and 
Smal\o\  in \zitat{\RRS}. They determined sets of polynomials cutting each 
of the irreducible components out of its ambient affine space.  This 
turned out to be quite taxing, even in the restricted scenario addressed. 
One of the  goals of the present article is to support the assertion that,
with respect to the equational aspect (among others) of the components, 
the projective setting has a significant edge over the affine.

Roughly, the underlying idea is as follows:  Given any irreducible 
component $\C$ of $\grassbd$, there is a representation-theoretically 
defined affine subvariety $\C(\S)$ of $\C$, which typically has very low 
dimension relative to $\dim \C$, such that: $\bullet$ $\C(\S) \times 
\AA^m$ is isomorphic to a dense open subset of $\C$ for some 
well-understood integer $m$, whence all birational data of $\C$ can be 
gleaned from $\C(\S)$;  $\bullet$ the ``dehydrated" format $\C(\S)$ of 
$\C$ holds a complete complement of information on $\C$, in that  it 
permits construction of a set of homogeneous polynomials defining $\C$ as 
a subset of the projective space $\PP$ in which $\grassbd$ lives (this 
construction follows a rote procedure which relies entirely on $\C(\S)$ 
and on the corresponding {\it skeleton\/} $\S$;  see below);  $\bullet$ 
polynomials which cut $\C(\S)$ out of a ``small" affine space are readily 
available from quiver and relations of $A$ by way of a direct 
combinatorial bridge, thus allowing for geometric and 
representation-theoretic analysis of $\C$ from a quiver presentation of 
$A$.  Polynomials for $\C$ in the ambient projective space $\PP$ are 
``inflated'' versions of those for $\C(\S)$.

Here is a red thread through the reduction process $\C \rightsquigarrow 
\C(\S)$.  The radical layering of the modules in an irreducible component 
$\C$ of $\grassbd$ is generic, i.e., constant on a dense open subset of 
$\C$. (Recall that the radical layering of an $A$-module $M$ is the 
sequence of semisimple modules $\SS(M) = (J^lM/J^{l+1}M)_{0 \le l \le L}$, 
where $J$ is the Jacobson radical of $A$ and $L+1$ the Loewy length.)   If 
the sequence $\SS = (\SS_0, \dots, \SS_L)$ of semisimples in $\lamod$ is 
the generic radical layering of $\C$ and $\biggrassSS$ denotes the locally 
closed subvariety of $\grassbd$ consisting of the points (representing 
modules) with radical layering $\SS$, then $\C \cap \biggrassSS$ is an 
irreducible component of $\biggrassSS$, and $\C$ is the closure of this 
intersection in the projective space $\PP$.  However, not all irreducible 
components of $\biggrassSS$ close up to irreducible components in 
$\grassbd$ --  such closures may fail to be maximal irreducible.   On the 
other hand, for any sequence $\SS$, polynomials that cut out the closures 
in $\PP$ of the irreducible components of the varieties $\biggrassSS$ are 
available from quiver and relations of $A$ via a simple algorithm.

To continue our outline, we now suppose that $\frakD$ is an irreducible 
component of $\biggrassSS$.  The key to further reduction consists of 
accessible affine subvarieties of $\frakD$;  see \zitat{\hierarchies} and 
\zitat{\BHT}.  They are defined in terms of distinguished path bases $\S$ 
of $A$-modules which are dubbed skeleta. The skeleta  --  finitely many 
for fixed $\SS$  --  are completely determined by the Gabriel quiver of 
$A$.   The corresponding varieties $\biggrassS$, each consisting of the 
points in $\biggrassSS$ that represent modules with skeleton $\S$, form an 
open affine cover of $\biggrassSS$.   Consequently, $\frakD$ coincides 
with the closure in $\biggrassSS$ of  any nontrivial intersection $\frakD 
\cap \biggrassS$.

The latter intersection can be whittled down further without any loss of 
birational information.
Namely, suppose that  $\S$ is a skeleton with $ \frakD \cap \biggrassS \ne 
\varnothing$.  This intersection is a direct product $\frakD(\S) \times 
\AA^{m(\SS)}$, where $m(\SS)$ is typically ``large", and $\frakD(\S)$ is 
the intersection of $\frakD$ with a  ``small" affine subvariety $\grassS$ 
of  the ambient projective space $\PP$;  in fact, $\frakD(\S)$ is an 
irreducible component of $\grassS$.  Polynomials for $\grassS$ are 
available from quiver and relations of $A$ almost at a glance, whence 
polynomials for the irreducible components of $\grassS$ can be computed 
(e.g., by using the Macaulay2 computer package).    In Section 6, we show 
how these polynomials lead to homogeneous equations for the closure of 
$\frakD$ in $\PP$.
On the side, we mention that polynomials for $\frakD(\S)$ also yield a 
minimal projective presentation of a ``generic module" for $\frakD$; see 
\zitat{\BHT, Section 4}.

Furthermore, we extend the available transfer of geometric information 
among the projective and affine settings, established in \zitat{\GeomIV}; 
see Diagram 2.1 for a table of correspondences.  Namely, we prove that 
smoothness and normality carry over from the irreducible components of the 
small varieties $\grassS$ to those of the varieties $\biggrassSS \subseteq 
\grassbd$ and $\repSS \subseteq \replabd$. Moreover, unirationality of the 
irreducible components of the $\grassS$ is inherited by those of 
$\grassbd$ and $\replabd$; see Theorems 3.2, 3.3.  Unirationality, in 
particular, entails representation-theoretic benefits:  On one hand, the 
mentioned presentation of generic modules becomes more transparent; on the 
other, the following result of Koll\'ar \zitat{\degen, Proposition 3.6} 
provides a tool towards determining non-generic representations in a given 
irreducible component from a generic one:   Namely, if $V$ is a 
unirational projective variety, then any finite set of points in $V$ 
belongs to some curve  $\PP^1 \rightarrow V$. Among the algebras for which 
the above properties are guaranteed are the truncated path algebras 
\zitat{\BHT} and the self-injective algebras with $J^3 = 0$ (see 
Proposition 5.2 and Theorem 5.4 below); the latter include the Carlson 
algebra.

A full list of polynomials pinning down an irreducible component in either 
$\replabd$ or $\grassbd$ is typically overwhelming.  The usefulness of 
such a list appears to be restricted to situations where the polynomials 
can be qualitatively understood sufficiently well to permit analysis of 
non-birational features, such as normality.   On the basis of a 
preliminary investigation,  we conjecture that all irreducible components 
of the varieties $\replabd$ and $\grassbd$ are normal in case $A$ is 
either a truncated path algebra or a self-injective algebra with $J^3 = 
0$.   (See Section 5 and the end of Section 6 for plausibility.)

To return to the overarching problem of understanding the irreducible 
components of $\replabd$ and $\grassbd$: The remaining (major) difficulty 
lies in the task of singling them out  among the closures of the 
irreducible components of the respective subvarieties $\repSS$ or 
$\biggrassSS$.   The above methods provide a framework for a ``generic 
representation theory" addressing representations with fixed radical 
layering $\SS$.  For fixed dimension vector $\bd$, they merely produce 
finite sets of  closed irreducible subvarieties of $\grassbd$ and 
$\replabd$, among which the irreducible components of the 
latter varieties are known to occur.

\definition{Outline of the article}  In Section 2, we give an overview of 
the relevant paramet\-rizing varieties and their interconnections; in 
particular, skeleta are defined, and  affine coordinates for the small 
varieties $\grassS$ parametrizing the representations with fixed skeleton 
$\S$ are described.  In Section 3, we establish the theorems ensuring 
transfer of smoothness, normality, and unirationality from irreducible 
components of the $\grassS$ to irreducible components of the larger 
parametrizing varieties.  In Section 4, we recall from \zitat{\BHT} and 
\zitat{\hierarchies} how to obtain polynomials  for the irreducible 
components of the pivotal varieties $\grassS$.   We list the generic 
skeleta for the  Carlson algebra in Section 5.  In the same section, we 
bring the criteria of Section 3 to bear:  For the larger class of 
self-injective algebras $A$ with $J^3 = 0$, we derive unirationality of 
the  components of the varieties $\replabd$, as well as smoothness of the 
components of the $\repSS$.  In Section 6, a general procedure is 
described for finding polynomial equations for  the closures of the 
irreducible components of the $\biggrassSS$ in the projective space $\PP$ 
in which $\grassbd$ is located. The method is illustrated by a small 
parametrizing variety over the Carlson algebra.   In particular, this 
procedure makes it evident that moving from the scenario of $\grassS$ to 
that of the projective space containing $\grassbd$ simply amounts to a 
change of bookkeeping.
\enddefinition

\definition{Conventions} 
Throughout, $K$ is an algebraically closed field, $A$ a basic finite dimensional $K$-algebra with Jacobson radical $J$, and $L+1$ the Loewy length of $A$. We do not lose generality in assuming that $A= KQ/I$, where $Q$ is a quiver and $I \subseteq KQ$ an admissible ideal in the path algebra. Our convention for multiplying paths $p$, $q$ in $KQ$ is to write $pq$ for ``$p$ after  $q$". The vertices $e_1,\dots,e_n$ for $Q$, abbreviated to $1,\dots,n$ in graphs, will be identified with a full set of primitive idempotents in either $KQ$ or $A$.

A {\it semisimple sequence with top $T$ and dimension vector $\bd$\/} is a sequence $\SS = (\SS_0, \dots, \SS_L)$ of semisimple objects in $\lamod$ such that $\SS_0 = T$ and $\underbardim \bigoplus_{0 \le l \le L} \SS_l = \bd$.  The {\it top\/} of a left $A$-module $M$ is $T(M) := M/JM$, and thus the {\it radical layering\/} of $M$, 
$$\SS(M) := (J^lM/J^{l+1}M)_{0 \le l \le L} \,,$$ 
is an example of  a semisimple sequence with top $T(M)$ and dimension vector $\underbardim M$. In dealing with tops and semisimple sequences, we systematically identify isomorphic semisimple modules.

A {\it top element\/} of $M$ is any element $x\in M\setminus JM$ such that $e_ix=x$ for some $i\le n$; in this situation, we say that $x$ is {\it normed by $e_i$\/}. Top elements $x_1,\dots,x_t \in M$ form a {\it full sequence of top elements of $M$\/} in case the residue classes $x_i+JM$ form a basis for $M/JM$.

All of the parametrizing varieties in the projective setting are subvarieties of Grassmannians of $K$-subspaces of the following projective $A$-modules $P$ and $\bP$.
Given a semisimple module $T$ of dimension $t$, we fix a projective cover $P := \bigoplus_{1\le r\le t} Az_r$ of $T$, together with a {\it distinguished\/} ({\it full\/}) {\it sequence\/} $z_1,\dots,z_t$ of top elements of $P$. By $\Phat$, we denote the corresponding projective $KQ$-module
$$\Phat := \bigoplus_{1\le r\le t} KQz_r \,.$$
A {\it path in $P$\/} (resp., {\it in $\Phat$\/}) is any nonzero element $pz_r \in P$ (resp., $pz_r\in \Phat$), where $p$ is a path in $Q$; in particular, the initial vertex of $p$ coincides with the vertex that norms $z_r$, and $p\in KQ\setminus I$ in case $pz_r\in P$.  In either case, it is unambiguous to refer to the starting and end points of $p$ as {\it starting and end points of the path $pz_r \in P$\/}.  Moreover, for a path $pz_r$ in $\Phat$, we may define $\len(pz_r) := \len(p)$.  It is the path length grading of $KQ$ that prompts us to shift back and forth between $A$- and $KQ$-modules.  We call $qz_r$ an {\it initial subpath\/} of the path $pz_r$ in $\Phat$ in case $p=uq$ for some path $u$ of length $\ge0$ in $KQ$.

We use the same notation for a path $pz_r$ in $\Phat$ as for the corresponding element of $P$. Note that, given a path $pz_r$ in $\Phat$, the canonical image $pz_r\in P$ is a path in $P$ (i.e., is nonzero) only when $p\notin I$.

Analogously, for a dimension vector $\bd= (d_1,\dots,d_n)$, we let $\bP := \bigoplus_{1 \le r \le d} A z_r$ be a fixed projective cover of $\bigoplus_{1 \le i \le n} S_i^{d_i}$ with distinguished sequence $z_1, \dots, z_d$ of top elements, and we denote by $e(r)$ the vertex of $Q$ which norms $z_r$.  Following the above pattern, we moreover introduce the $KQ$-module $\bPhat$.   In particular, this provides us with canonical embeddings $P \subseteq \bP$ and $\Phat \subseteq \bPhat$.   The concepts of a {\it path in\/} $\bP$ or $\bPhat$ are carried along as well. For any permutation $\tau$ of $\{1,\dots,d\}$ such that $e(\tau(r))= e(r)$ for all $r$, the $KQ$-auto\-morph\-ism of $\bPhat$ or $\bP$ defined by $z_r \mapsto z_{\tau(r)}$  is called a {\it permutation automorphism\/}.  Notationally, we do not distinguish between the actions on $\bP$ and $\bPhat$.
\enddefinition

\head 2.  Prerequisites concerning the affine and projective module varieties \endhead

\definition{A.  Overview}  We start with an overview, Diagram 2.1, of all the relevant varieties (slightly updating the notation of \zitat{\hierarchies}), and follow with definitions. Then we state how information is transferred. In particular, we explain how the irreducible components of the displayed affine and projective varieties arise as successive closures of irreducible components in the smallest ones.
\goodbreak \midinsert
$$\xymatrixrowsep{0.2pc}\xymatrixcolsep{3pc}
\xymatrix{
\txt{``big'' scenario}  &&\txt{``small'' scenario}  \\
\boxed{\grassbd} \dropdown{\txt{(projective)}} \ar@{<->}[r] &\boxed{\replabd}
\dropdown{\txt{(affine)}} \\  \\
\smvsubseteq &\smvsubseteq \\
\boxed{\biggrasstbd} \dropdown{\txt{(quasi-projective)}} \ar@{<->}[r]
&\boxed{\reptbd} \dropdown{\txt{(quasi-affine)}}
\ar@{<->}[r] &\boxed{\grasstbd} \dropdown{\txt{(projective)}} \\  \\
\smvsubseteq & \smvsubseteq & \smvsubseteq \\
\boxed{\biggrassSS} \dropdown{\txt{(quasi-projective)}} \ar@{<->}[r]
&\boxed{\repSS} \dropdown{\txt{(quasi-affine)}}
\ar@{<->}[r] &\boxed{\grassSS} \dropdown{\txt{(quasi-projective)}} \\  \\
\smvsubseteq && \smvsubseteq \\
\boxed{\biggrassS} \dropdown{\txt{(affine)}} \ar@{<->}[rr] &&\boxed{\grassS}
\dropdown{\txt{(affine)}}
 }$$
\medskip
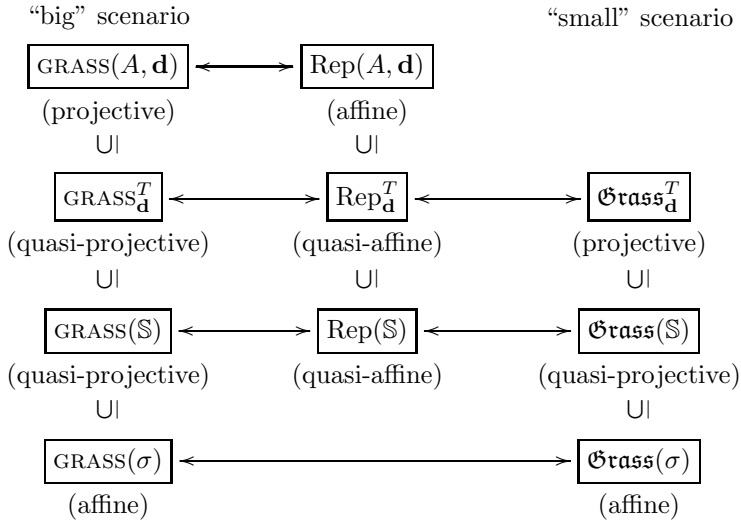
\centerline{Diagram 2.1}
\endinsert

The varieties in the left-most column of Diagram 2.1 are subvarieties of the Grass\-mann\-ian of $(\dim\bP-d)$-dimensional $K$-subspaces of $\bP$, those in the right-most columns are subvarieties of the Grassmannian $\Gr(\dim P-d, P)$ with $P$ and $\bP$ as above; here $d = |\bd|$.  We introduce the  displayed varieties from top to bottom, moving from right to left in each row. Let  $T$  be a semisimple $A$-module with dimension vector $\le\bd$, and $\SS$ a semisimple sequence with top $T$ and dimension vector $\bd$.

\roster
\item $\replabd$: This is the classical affine variety parametrizing the left $A$-mod\-ules with dimension vector $\bd$, namely

\noindent\qquad$\bigl\{ (x_\alpha)_{\alpha \in Q_1} \in \prod_{\alpha \in Q_1} \Hom_K \bigl(K^{d_{\text{start}(\alpha)}},\, K^{d_{\text{end}(\alpha)}}\bigr) \mid \text{the\ } x_{\alpha}$

\rightline{satisfy all relations in $I \bigr\}$,\qquad}

\noindent where $Q_1$ is the set of arrows of the quiver $Q$.
 As usual, we endow $\replabd$ with the conjugation action of $\GLbd := \GL_{d_1}(K) \times \cdots\times \GL_{d_n}(K)$.
\item $\grassbd$: Set $a := \dim\bP-d$. Then $\grassbd$ is the closed subvariety of the Grassmannian $\Gr(a,\bP) \subseteq \PP(\la^a\bP)$ consisting of those $K$-subspaces $C\subseteq \bP$ which are $A$-submodules of $\bP$ such that $\underbardim(\bP/C)=\bd$. In particular, $\grassbd$ is a projective variety. Note, moreover:  The linear algebraic group $\aut_{A}(\bP)$ acts mor\-phically on  $\grassbd$, and the orbits of this action are in 1--1 correspondence with the isomorphism classes of left $A$-modules having dimension vector $\bd$.
\item $\grasstbd$: This denotes the closed subvariety of the classical Grassmannian $\Gr(\dim P-d, P)$ (which is typically far smaller than $\Gr(a,\bP)$) consisting of those $K$-subspaces $C\subseteq P$ which are $A$-submodules of $P$ with the additional properties that $\underbardim(P/C)= \bd$ and $P/C$ has top $T$.  Observe that the latter condition implies that all points of $\grasstbd$ are submodules of $JP$. Clearly, $\grasstbd$ carries a morphic action of the smaller automorphism group $\autlap$.  This time, the orbits of $\grasstbd$ under the action are in 1--1 correspondence with the isomorphism classes of left $A$-modules with dimension vector $\bd$ and top $T$.  In light of the embedding $P \subseteq \bP$, we note that $\grasstbd$ embeds into $\grassbd$, via $C \mapsto C \oplus \bigoplus_{t +1 \le r \le d} A z_r$.  This embedding makes $\grasstbd$ a closed subvariety of $\grassbd$ which, however, fails to be closed under the action of $\aut_A(\bP)$ in general. 
\item $\reptbd$: This is the locally closed subvariety of $\replabd$ consisting of the points corresponding to modules with top $T$. It is clearly closed under the $\GLbd$-action.
\item $\biggrasstbd$ is the locally closed subvariety of $\grassbd$ consisting of all points in the latter variety that correspond  to modules with top $T$. Clearly,  the embedding $\grasstd \hookrightarrow \grassbd$ under (3) makes $\grasstd$ a closed subvariety of $\biggrasstbd$.  We observe stability of $\biggrasstbd$ under the $\aut_{A}(\bP)$-action; in fact $\biggrasstbd$ is the closure of $\grasstd$ under the action of the big automorphism group.  
\item $\grassSS$, $\repSS$, $\biggrassSS$: Each of these is the locally closed subvariety of the variety shown above it that consists of the points parametrizing the modules with radical layering $\SS$. Evidently, in each case, the mentioned subvariety is closed under the pertinent group action, that of $\autlap$, $\GLbd$, and $\aut_{A}(\bP)$, respectively.  Observe that $\biggrassSS$ is just the closure of  $\grassSS$ under  the $\aut_{A}(\bP)$-action.  The varieties $\grassSS$ and $\biggrassSS$ are very similar geometrically in that they have open affine covers, the patches of which differ only by a direct factor $\AA^m$.   (However, compared with $\grassSS$, the higher dimension of $\biggrassSS$ leads to a plethora of additional polynomial equations as one passes to the closure in the ambient projective space $\PP = \PP(\la^a\bP)$.  These additional equations are difficult to analyze in terms of their geometric implications,  such as appearance of  new singularities,   in general.)
\item $\grassS$, $\biggrassS$: The definitions of these varieties hinge on the concept of a skeleton (see below), as introduced and studied in the present generality in \zitat{\hierarchies} and \zitat{\BHT}.  Namely, $\grassS$ and $\biggrassS$ are the respective open subvarieties of $\grassSS$ and $\biggrassSS$ consisting of the points that represent modules with skeleton $\S$.
\endroster
\enddefinition

\definition{B. Skeleta} Suppose the semisimple sequence $\SS$ with top $T$ and dimension vector $\bd$ has the form $\SS = (\SS_0, \dots, \SS_L)$. The terms of the sequence will be conveyed in the format $\SS_l = \bigoplus_{1 \le i \le n} S_i^{m(l, i)}$. Moreover, we let $P$, $\bP$, $\Phat$, $\bPhat$ and $z_1, \dots, z_d$ be as at the end of Section 1.  Recall that our embedding  $\Phat \subseteq \bPhat$ is such that the first $t$ members of the distinguished sequence of top elements of $\bPhat$ are the distinguished top elements of $\Phat$.
\medskip

{\bf {(i)}} An ({\it abstract\/}) {\it skeleton in $\Phat$\/} (resp., {\it in $\bPhat$\/}) {\it with  radical layering $\SS$\/}  is a set $\S$ of paths in $\Phat$ (resp., in $\bPhat$) of lengths at most $L$  with the following two properties:  

$\bullet$  For each $l \in \{0 , \dots ,L\}$ and each $i \in \{1, \dots, n\}$, the set $\S$ contains precisely $m(l,i)$ paths of length $l$ ending in $e_i$.  (In particular, $\S$
contains precisely $t$ top elements $z_{i_1}, \dots z_{i_t}$.  In case $\sigma \subseteq \Phat$, these $t$ top elements are $z_1, \dots, z_t$.)

$\bullet$ $\S$ is closed under initial subpaths, that is, whenever
$p_2 p_1 z_r \in \S$, then $p_1 z_r \in \S$. 
\smallskip

We write $\SS(\S) = \SS$ and $\underbardim\S= \underbardim\SS= \bd$ in this situation. 
\smallskip

\noindent Alternatively, we view a skeleton
$\S$ as a forest of $t$ tree graphs, where the $r$-th tree displays, as edge paths starting at the root, the paths
in $\S$ starting in the top element $z_{i_r}$.  For examples of graphs of skeleta, we refer to Sections 4,5 below and to \zitat{\hierarchies}.
\medskip    

{\bf {(ii)}}  Let $\S$ be a skeleton in $\Phat$ (resp., in $\bPhat$) and $M$ an
$A$-module with $\SS(\S) = \SS(M)$.  Moreover, we denote by  $i_1<\cdots <i_t$ the indices $i\in \{1,\dots,d\}$ such that $z_i\in \S$.

We call $\S$ a {\it skeleton
of\/} $M$ if there exists a full sequence of top elements $x_1, \dots,
x_t$ of $M$, together with a $KQ$-epimorphism $f: \Phat \rightarrow
M$ (resp., $f: \bPhat \rightarrow M$) such that $f(z_{i_r}) = x_r$ for all $r$ and the set $f(\S)$ is a $K$-basis for $M$.    In this
situation, we also say that $\S$ is a {\it skeleton of $M$ relative to
the sequence $x_1, \dots, x_t$\/}.
In the special case where $M = P/C$ for some point $C \in \grass(\SS(\S))$ (resp., $M= \bP/C$ for some $C\in \GRASS(\SS(\S))$), we say that $\S$ is a {\it distinguished skeleton of $M$\/} if $\S$ is a skeleton of $M$ relative to the distinguished sequence of top elements $z_{i_1} + C, \dots, z_{i_t} + C$.
\medskip

{\bf {(iii)}} Again, let $a := \dim \bP - d$. For a skeleton $\S$ in $\Phat$ (resp., in $\bPhat$), we define $\grassS$ (resp.,  $\biggrassS$ and $\SchubertS$), as follows:
$$\align
\grassS  &:= \{ C \in \grass(\SS(\S)) \mid \S \text{\ is a distinguished
skeleton of\ } P/ C \}  \\
\biggrassS &:= \{ C \in \GRASS(\SS(\S)) \mid \S \text{\ is a distinguished
skeleton of\ }  \bP/C \}  \\
\SchubertS &:= \{C \in \Gr(a,\bP) \mid \bP = \bigl(  \bigoplus_{b \in \S} K b \bigr) \oplus C\}.
\endalign$$
\enddefinition 

\definition{Comments} 
{\bf Re (i).}  In the definition, we place $\S$ into $\Phat$ or $\bPhat$ rather than $P$ or $\bP$, in order to have an unambiguous notion of path length.  When there is no risk of ambiguity, we typically identify $\S$ with its canonical image in $P$ or $\bP$.  Whenever $\biggrassS \ne \varnothing$, this identification is harmless, since the image of $\S$ in $P$ (resp., $\bP$) is still $K$-linearly independent in that case.

Any skeleton $\S$ in $\bPhat$ can be mapped to a skeleton $\S'= g.\S$ in $\Phat$ by a permutation automorphism $g$ of $\bPhat$.
\smallskip  

{\bf Re (ii).} Since we insist that $\S$ and $M$ have the same radical layering, the requirement that the set $\{ p x_r \mid p z_{i_r} \in \S\}$ be a basis for $M$ amounts to the following:  For each $l \in \{0, \dots, L\}$, the set
$$\{p x_r \mid p z_{i_r} \in \S,\, \len(pz_{i_r}) = l\} $$  
induces a $K$-basis for the radical layer $J^lM / J^{l+1} M$. Thus the skeleta of $M$ are $K$-bases which
are tied to the $KQ$-structure of $M$.  Clearly, the set of skeleta of any $M \in \lamod$ is nonempty and finite; this set is an isomorphism invariant of $M$.
\smallskip

{\bf Re (iii).} If $\S$, $\S'$ are skeleta in $\bPhat$ and $\S'= g.\S$ for some permutation automorphism $g$ of $\bPhat$, then $g$ (viewed as an automorphism of $\bP$) maps $\biggrassS$ isomorphically onto $\GRASS(\S')$. This allows us to concentrate on skeleta in $\Phat$ when we explore $\aut_A(\bP)$-stable subsets of $\grassbd$.

Observe that  $\biggrassS = \SchubertS \cap \biggrassSS$. Analogously, $\grassS$ is the intersection of $\grassSS$ with the correspondingly defined smaller Schubert cell of $\PP(\la^{\dim P  -d} P)$ determined by $\S$.  While $\SchubertS$ is clearly an open subvariety of $\PP(\la^a\, \bP)$, the intersection $\biggrassS$ need not be open in $\grassbd$.
However, in view of \zitat{\hierarchies, Section 3}, we know:

$\bullet$ The varieties $\grassS$, where $\S$ traces the skeleta in $\Phat$ with radical layering $\SS$, form an  {\it open\/} affine cover of $\grassSS$.  Each $\grassS$ is stable under the action of the unipotent radical of $\autlap$, but not under the action of $\autlap$, in general.

$\bullet$  The varieties $\biggrassS$, as $\S$ traces the skeleta in $\bPhat$ with radical layering $\SS$, form an  {\it open\/} affine cover of $\biggrassSS$.  Moreover, $\biggrassS$ is stable under the action of the unipotent radical of $\aut_A(\bP)$.  In case $\S \subseteq \Phat$, we have $\biggrassS\cong \grassS \times \AA^{m}$, where $m = \sum_{1 \le i \le n} d_i(d_i - t_i)$ and $(t_1,\dots,t_n)= \underbardim T$. In view of these isomorphisms, the irreducible components of $\biggrassS$ are isomorphic to the varieties $\C(\S)\times \AA^{m}$ as $\C(\S)$ traces the irreducible components of $\grassS$.

In particular, the above facts guarantee the set of skeleta of a module to be generically constant on the irreducible components of all of the considered parametrizing varieties.
\enddefinition

\definition{C. Transfer of information among the displayed varieties}
First we recall the horizontal connections among the irreducible components of the varieties in Diagram 2.1. 
\enddefinition  

\proclaim{Theorem 2.1} {\rm\zitat{\GeomIV, Proposition C}}.

{\bf (I)}  Consider the one-to-one correspondence between the orbits of $\grasstbd$ and $\reptbd $ which
assigns to any orbit $\autlap.C \subseteq \grasstbd$ the orbit
$\GLbd.x \subseteq \reptbd $  representing the same
$A$-module up to isomorphism.  This correspondence extends to an inclusion-pre\-serv\-ing bijection
$$\Phi: \{ \autlap\text{-stable subsets of\ } \grasstbd \} \rightarrow
\{\GLbd\text{-stable subsets of\ } \reptbd \}$$  
which preserves and reflects
openness, closures, connectedness, irreducibility, and types of
singularities.
\smallskip

{\bf (II)}  The one-to-one correspondence which
assigns to any orbit $\aut_A(\bP).C$ of $\grassbd$ the orbit
$\GLbd.x$ of $\replabd$  representing the same isomorphism class of $A$-modules extends to an inclusion-preserving bijection
$$\multline
\Phi: \{ \autlabp\text{-stable subsets of\ } \grassbd \} \rightarrow \\
\{\GLbd\text{-stable subsets of\ } \replabd \}
\endmultline$$  
which preserves and reflects
openness, closures, connectedness, irreducibility, and types of
singularities.
\qed
\endproclaim

Theorem 2.1 applies, in particular, to the irreducible components of the varieties displayed in the top three rows of Diagram 2.1. Indeed, due to connectedness of the groups $\autlap$, $\autlabp$ and $\GLbd$, the irreducible components of these varieties are stable under the operation of the pertinent acting groups.  (To see this, note that, given an irreducible component $\C$ of $\grassSS$ for instance, the image of the canonical morphism $\autlap \times \C \rightarrow \grassSS$ is irreducible.)  Consequently, we have geometrically useful one-to-one correspondences among the sets of irreducible components along the horizontal double arrows in the top three rows. We apply these correspondences in Sections 3, 5, 6.

Recall that a ``type of singularity" is an equivalence class of pointed varieties $(X,x)$ under the following equivalence relation ({\it smooth equivalence\/}):  $(X,x) \sim (Y,y)$ precisely when there exists a pointed variety $(Z,z)$ together with smooth morphisms $Z \rightarrow X$ and $Z \rightarrow Y$ taking $z$ to $x$ and $y$, respectively \zitat{\Hess, (1.7)}.  By \zitat{\EGAfourtwo,  Sections 5,6}, many local properties are invariant under this equivalence relation. In particular, given $(X,x)\sim (Y,y)$, the point $x$ is smooth (resp., normal) if and only if $y$ is smooth (resp., normal) \zitat{ibid} (see also \zitat{\djong, Commutative Algebra, Lemmas 138.7, 139.3} for normality).

\definition{Vertical connections} 
Next, we address the (known) vertical connections, tracking the irreducible components from the bottom to the top of the columns in the diagram; the easy proofs can be found in \zitat{\BHT, Section 2}.   

 The irreducible components of $\grassSS$ (resp., $\biggrassSS$) are precisely the closures in $\grassSS$ (resp., $\biggrassSS$) of the irreducible components of the $\grassS$ (resp., $\biggrassS$). In other words: the problem of identifying the irreducible components of $\grassSS$ is equivalent to decomposing the varieties $\grassS$ into their irreducible components, where $\S$ traces the skeleta with radical layering $\SS$;  analogously for $\biggrassSS$.

The irreducible components of $\grasstbd$ are among the closures in $\grasstbd$ of the irreducible components of the $\grassSS$, but the latter closures may fail to be maximal irreducible in $\grasstd$. The same holds for irreducible components of $\biggrasstbd$ and those of $\biggrassSS$.

Finally, the irreducible components of $\grassbd$ are among the closures in $\grassbd$ of the irreducible components of the $\biggrasstbd$, with the same caveat as above.

While the varieties $\grassS$ do not have a useful counterpart in $\replabd$, all of the above connections among the irreducible components of the varieties in the Grassmannian setting carry over, mutatis mutandis, to those of the varieties $\replabd$ versus those of $\reptbd$ and $\repSS$; this is immediate from the remarks following Theorem 2.1.

The unresolved problem of classifying the irreducible components of the module varieties $\replabd$ and $\grassbd$ for a given algebra $A$ thus consists of finding the generic radical layerings, and the generic skeleta among those having generic radical layering;  finally, of identifying among the irreducible components of the pertinent $\grassS$ those which close up to maximal irreducible sets in $\grassbd$. 
\enddefinition

\definition{D. Affine coordinates for the varieties $\grassS$ and $\biggrassS$}
Once more, we let $\S \subseteq \Phat$ be an abstract skeleton with radical layering $\SS$, where $\SS$ is a semisimple sequence with top $T$ and dimension vector $\bd$.  

A {\it $\S$-critical\/} path in $\Phat$ is any path $b'$ of length at most $L$ in $\Phat \setminus \S$ with the property that every proper initial subpath belongs to $\S$. Since the distinguished top elements of $\Phat$ all lie in $\S$, the criticality condition means that $b' = \alpha pz_r \notin \S$, where $\alpha$ is an arrow and $pz_r \in \S$ has length $< L$.  Moreover, to any $\S$-critical path $b'$ we associate the following subset of $\S$:  
$$\S(b') := \{ b \in \S \mid \len(b) \ge \len(b') \text{\ and\ end}(b) = \text{end}(b') \}.$$
As before, we use the same notation for the incarnations of a path $pz_r$ in both $\Phat$ and $P$.

Finally, we let $N$ be the (disjoint) union of the sets  $\{b' \} \times \S(b')$, where $b'$ traces the $\S$-critical paths.  

To any point $C \in \grassS$,  we now assign the element $(\cbb)_{(b', b) \in N} \in \AA^N$ which is determined by the following equalities in the factor module $P/C$:  
$$\xalignat2
 b' + C &=  \sum_{b \in \S(b')} \cbb\, b+C  &&(\S\text{-critical paths\;} b').
\tag 2.1 \endxalignat $$

According to \zitat{\hierarchies, Theorem 3.12},   the above assignment defines an isomorphism from  $\grassS$ onto a closed subvariety of $\AA^N$.  Polynomial equations for the latter can be algorithmically obtained from quiver and relations of $A$ \zitat{\hierarchies, Subsection 3.14}; the algorithm is recalled in Section 4, as it is the crucial link in obtaining polynomial equations for the components of the variety $\biggrassSS$ in its ambient projective space. 

Given a skeleton $\S$ in $\bPhat$, there is a parallel affine coordinatization for the subvariety $\biggrassS$ of $\grassbd$. Let $z_{i_1},\dots,z_{i_t}$ be the paths of length $0$ in $\S$. We carry over the definition of a {\it $\S$-critical path\/} essentially verbatim, merely replacing $\Phat$ by $\bPhat$. Note, however, that $z_r$ is $\S$-critical for all $r\notin \{i_1,\dots,i_t\}$ in the large scenario. The sets $\S(b')$, for $\S$-critical paths $b'$, are defined as before. In particular, for $r\notin \{i_1,\dots,i_t\}$, the set $\S(z_r)$ consists of {\it all\/} paths in $\S$ which end in the vertex $e(r)$ that norms $z_r$. In complete analogy with the small scenario, one obtains an isomorphism from $\biggrassS$ to a closed subvariety of $\AA^{\bold{N}}$, where $\bold{N}$ is the enlarged (disjoint) union of the sets $\{b'\} \times \S(b')$ with $b'$ tracing the $\S$-critical paths. Embeddings of $\biggrassS$ into an affine space of dimension $\le |\bold N|$ and into a suitable Schubert cell are described in Section 6.
\enddefinition

\head 3.  Transfer of properties from the $\operatorname{\frak{Grass}}(\sigma)$ to the larger parametrizing varieties \endhead

The purpose of this section is to demonstrate how properties of the $\grassS$, the smallest and most accessible of the parametrizing varieties, carry over to the irreducible components of the varieties $\grassSS$, $\repSS$, $\grassbd$, $\replabd$.  We focus on unirationality, normality  and smoothness.  First applications will follow in Sections 5 and 6.

Recall that an irreducible variety over $K$ is {\it unirational\/} if its rational function field embeds into a purely transcendental extension of $K$.   In light of openness of each $\grassS$ in $\grassSS$, unirationality of any irreducible component $\C$ of $\grassSS$ entails unirationality  of any nonempty intersection $\C \cap \grassS$ and vice versa.   In particular, unirationality of all components of the $\grassS$ with $\SS(\S) = \SS$ guarantees unirationality of all components of $\grassSS$. This condition also yields unirationality of all components of $\biggrassSS$, since the components of the corresponding $\biggrassS$ have the form $\C\nomathbreak \times\nomathbreak \AA^m$ for components $\C$ of $\grassS$.  Similarly, unirationality of the components of the $\biggrassSS$, for all semisimple sequences $\SS$ with dimension vector $\bd$, implies unirationality of the components of $\grassbd$;  likewise for unirationality of the components of $\repSS$ versus $\replabd$. Hence, the only point in the upcoming theorem that needs additional justification is the transit from the Grassmannian to the classical affine scenario.

We first single out a useful lemma.

\proclaim{Lemma 3.1} Let $\S$ be a skeleton in $\Phat$ such that $\underbardim\S = \bd$ and $\grassS$ is nonempty. There is an isomorphism $\phi$ from $\grassS$ onto a closed subvariety of $\replabd$ such that for each $C\in \grassS$, the points $C$ and $\phi(C)$ parametrize isomorphic $A$-modules. 
\endproclaim

\demo{Proof} For $1 \le i \le n$, we first label a preferred basis for $K^{d_i}$ by those paths $b \in \S$ which end in the vertex $e_i$; there are precisely $d_i$ such paths.  

Now any point $C\in \grassS$ is identified with the unique point $c =  (\cbb) \in \AA^N$ for which the equations (2.1) hold,  where $N$ is defined as in Section 2.D. The $A$-module $P/C$ has dimension vector $\bd$, and each vector space $e_i(P/C)$ has a basis $\{b+C \mid b\in\S,\; \text{end}(b)= e_i \}$; we identify $e_i(P/C)$ with $K^{d_i}$ so that $b+C$ corresponds to $b$. We send the point $C$ to the point $\phi(C)= (x^c_{\alpha})_{\alpha \in Q_1}$ in $\replabd$ determined by the following property: for any arrow $\alpha: e_i \rightarrow e_j$, the map $x^c_\alpha \in \Hom_K \bigl( K^{d_i},\, K^{d_j}\bigr)$ is given, in terms of our chosen basis, by $b+C \mapsto \alpha b+C$, that is,
$$x^c_\alpha(b) = \cases \alpha b &(\text{if\;} \alpha b \in \S)\\  0 &(\text{if\;} \len(b)=L)\\  \sum_{s \in \S(\alpha b)} c_{\, \alpha b, s} \,  s &(\text{if\;} \alpha b \text{\;is\;} \S\text{-critical}). \endcases$$
We thus obtain a morphism $\phi: \grassS \rightarrow \replabd$. The stated properties are clear.
\qed\enddemo

\proclaim{Theorem 3.2. Unirationality}  

{\bf (I)}  Let $\SS$ be a semisimple sequence, $\C$ an irreducible component of $\grassSS$, and $\D$ the corresponding irreducible component of $\repSS$.  For any skeleton $\S$ in $\Phat$ such that $\C \cap \grassS \ne \varnothing$, we have:

{\bf (1)}  $\C$ is unirational if and only if $\C \cap \grassS$ is unirational.

{\bf (2)} If $\C$ is unirational, then  so is $\D$.
\smallskip

{\bf (II)}   Now suppose that $\C$ is an irreducible component of  $\grassbd$ and that $\D$ is the corresponding irreducible component of $\replabd$.  Let $\SS$ be the generic radical layering of $\C$, $\D$, and $\S$ any skeleton in $\Phat$ with radical layering $\SS$ such that $\C \cap \grassS \ne \varnothing$.  Then:

{\bf (1)}  $\C$ is unirational if and only if $\C \cap \grassS$ is unirational.

{\bf (2)} If $\C$ is unirational, then  so is $\D$.
\smallskip

{\bf (III)} Suppose that the irreducible components of $\grassS$ are unirational for all skel\-eta $\S$ in $\Phat$ with $\underbardim \S= \bd$.  Then all irreducible components of $\replabd$ and $\grassbd$ are unirational, as are the components of the varieties $\repSS$ and $\biggrassSS$ with $\underline{\dim}\, \SS = \bd$.
\endproclaim

\definition{Remark}  We point out that, in part {\bf (II)}, the choice of $\,\S \subseteq \Phat$ satisfying $\C \cap \grassS \ne \varnothing\,$ is non-restrictive:  Indeed, recall that, for any $\aut_A(\bPhat)$-stable subvariety $V$ of $\grassbd$ and any skeleton $\S \subseteq \bPhat$ with $V \cap \biggrassS \ne \varnothing$, there exists a permutation automorphism $g$ of $\bPhat$ such that $g. \S \subseteq \Phat$.  Clearly, this implies  $V \cap \grass(g. \S) \ne \varnothing$; see Section 2.B. \enddefinition

\demo{Proof}  {\bf (I)} As noted above, only {\bf(2)} needs proof. Assume that $\C(\S):= \C \cap \grassS$ is nonempty and unirational, and let $\phi: \grassS \rightarrow \replabd$ be the morphism of Lemma 3.1. Then $\D(\S): = \phi(\C(\S))$ is a closed subvariety of $\replabd$ which is isomorphic to $\C(\S)$.  Consequently, $\D(\S)$ is in turn unirational.     
In order to deduce unirationality of $\D$, let $V$ be the closure of $\C(\S)$ in $\grassSS$ under the
$\autlap$-action, and $W$ the closure of $\D(\S)$ under the $\GLbd$-action.  Then $V$ and $W$ are paired by the correspondence of Theorem 2.1.  Since $\C(\S)$ is open and dense in $\C$, we find the $\autlap$-stable subvariety 
$$V = \bigcup_{f \in \autlap} f.\C(\S)$$ 
of $\C$ to be open and dense as well.  From Theorem 2.1, we thus deduce that the  $\GLbd$-stable subvariety $W$ of $\D$ is in turn open and
dense.  The construction of $W$ as the closure
of $\D(\S)$ under the $\GLbd$-action therefore provides us with a dominant
morphism $\GLbd \times\, \D(\S) \rightarrow \D$.  
This implies that  $\D$ is unirational as postulated. 
\smallskip

The claims under {\bf (II)} now follow in light of the remarks preceding the lemma, and those under {\bf (III)} are immediate from {\bf (I)} and {\bf (II)}. \qed
\enddemo 

Of course, all conclusions of Theorem 3.2 remain true as long as there are enough unirational patches $\biggrassS$ to cover $\grassbd$. The same is true for the conclusions of Theorem 3.3, with unirationality replaced by smoothness or normality. 

As for normality and smoothness, these properties transfer pointwise via Theorem 2.1;  cf\. the comments following that theorem.  In view of the fact that, for any semisimple sequence $\SS$, the variety $\grassSS$ is covered by open subvarieties $\grassS$, where $\S$ traces the skeleta with $\SS(\S) = \SS$, we glean the following criteria.  They will come to bear in Sections 5 and 6.    

\proclaim{Theorem 3.3. Normality and smoothness}
  
{\bf (I)}  Let $\SS$ be a semisimple sequence, $\C$ an irreducible component of $\grassSS$, and $\D$ the corresponding irreducible component of $\repSS$.

Then  $\D$ is normal {\rm(}resp., smooth{\rm)} if and only if $\C$ is normal {\rm(}resp., smooth{\rm)}, if and only if all nonempty intersections $\C \cap \grassS$, for skeleta $\S$ in $\Phat$ with radical layering $\SS$, are normal {\rm(}resp., smooth{\rm)}.
\smallskip

{\bf (II)}   Now suppose that $\C$ is an irreducible component of  $\grassbd$ and $\D$ the corresponding irreducible component of $\replabd$.  

Then $\D$ is normal {\rm(}resp., smooth{\rm)} precisely when $\C$ is normal {\rm(}resp., smooth{\rm)}.   Moreover, if $\SS$ is the generic radical layering of the modules in $\C$, $\D$, the following are equivalent:

{\bf (1)}  All non-normal {\rm(}resp., non-smooth{\rm)} points of $\D$ belong to $\D \setminus \repSS$.

{\bf (2)}  All non-normal {\rm(}resp., non-smooth{\rm)} points of $\C$ belong to $\C \setminus \biggrassSS$.

{\bf (3)}  All nonempty intersections $\C \cap \biggrassS$, for skeleta $\S$ in $\Phat$ with radical layering $\SS$, are normal {\rm(}resp., smooth{\rm)}. 
\endproclaim

\demo{Proof} The equivalences in part {\bf (I)}, and the first equivalence of part {\bf (II)}, follow directly from Theorem 2.1 and the preceding discussion, as does the equivalence of {\bf (1)} and {\bf (2)}. Similarly, {\bf (2)} is equivalent to normality (resp., smoothness) of all nonempty intersections $\C \cap \GRASS(g.\S)$, for all skeleta $\S$ in $\Phat$ and all permutation automorphisms $g$ of $\bPhat$. Since $\C \cap \GRASS(g.\S) \cong \C\cap \biggrassS$ for any such $g$, we conclude that condition {\bf (2)} is equivalent to {\bf (3)}.
\qed\enddemo

On the other hand, if $\C$ is an irreducible component of $\grassbd$ as in part {\bf (II)} of the theorem, then smoothness of all nonempty intersections  $\C  \cap \grassS$ need not imply smoothness of $\C$.  This phenomenon will be exemplified in Section 6.

\head 4.  The varieties $\operatorname{\frak{Grass}}(\sigma)$ \endhead

As pointed out earlier, the varieties $\grassS$, where $\S$ traces the skeleta with dimension vector $\bd$, encode geometric information about $\grassbd$  --  and about $\replabd$ by extension  --  in the most compressed format.  An algorithm for deriving polynomials for the $\grassS$ in an appropriately small affine space is available in \zitat{\codes}.  We begin this section with a summary of the underpinnings proved in \zitat{\hierarchies}, in order to exhibit the direct connection between $\grassS$ and the relations of $A$. We end by illustrating the theory for the Carlson algebra.

Fix a semisimple sequence $\SS$ with top $T$ and dimension vector $\bd$, the distinguished projective cover $P = \bigoplus_{1 \le r \le t} A z_r$ of $T$, and a skeleton $\S$ in $\Phat$ with radical layering $\SS$.  We refer to Sections 1 and 2.D for further notation and terminology, in particular for the isomorphism $C \mapsto (\cbb)_{(b', b) \in N}$ from $\grassS$ onto a closed subvariety of $\AA^N$.   Here
$$N := \{ (b', b) \mid b' \text{\; a\;} \S\text{-critical path in\;} \Phat,\; b\in \S(b') \},$$
and the $\cbb$ are the unique scalars such that
$$C = \sum_{b' \, \S \text{-critical}} A \biggl (b' \, - \, \sum_{b \in \S(b')} \cbb\, b \biggr).$$
We also write $\grassS$ for the corresponding isomorphic copy contained in $\AA^N$.

The announced algorithm for obtaining polynomials defining $\grassS$ in $\AA^N$ is simply a reformulation of the fact that the following left module $F$ over the polynomial ring 
 $KQ[X] = KQ[ X_{b', b} \mid (b',b) \in N ]$
 is free over the smaller polynomial ring 
 $$K[X] = K[ X_{b', b} \mid (b',b) \in N ];$$ 
see \zitat{\hierarchies, Lemma 3.13}:  Namely, the factor module
 $$F = \bigl(KQ[X] \otimes _{KQ} \Phat\bigr)\, / \, \U,$$
 where $\U \subseteq KQ[X] \otimes _{KQ} \Phat$ is the left $KQ[X]$-submodule generated by the differences 
 $$b' - \sum_{b \in \S(b')} X_{b', b}\, b \qquad\quad (\S\text{-critical\;} b'\in \Phat),$$
together with all paths of length $L+1$ in $\Phat$.  A basis for $F$ over $K[X]$ is $\S$.  
 
 Let $R\subseteq \bigcup_{1\le k,l\le n} e_kIe_l \subseteq KQ$ be a finite generating set for the ideal $I\subseteq KQ$, viewed as a left ideal.  Recalling that the idempotent norming the top element $z_m$ of $\Phat$ is denoted by $e(m)$, we set
$$\align
\Rhat \,:=\, \{ \rho z_m \mid \rho\in R\,e(m),\; m=1,\dots,t\} \,&\subseteq\, \Phat  \\
 &\subseteq\, KQ[X] \otimes_{KQ} \Phat \,=\, \bigoplus_{1\le r\le t} KQ[X]z_r \, .
 \endalign$$
It will be convenient to enumerate the elements of $\S$, say $\S= \{b_1,\dots, b_d\}$.

Writing $\seq$ for congruence modulo $\U$ on $KQ[X] \otimes_{KQ} \Phat$, we obtain expansions
$$\rho z_m \seq \sum_{1\le l\le d} \tau_{\rho z_m,l}\, b_l \quad \text{with} \quad \tau_{\rho z_m,l} \in K[X]\,$$
where the polynomials $\tau_{\rho z_m,l}$ are uniquely determined by $\rho$ and $m$. The incarnation of $\grassS$ in $\AA^N$ is the vanishing locus of the resulting polynomials
$$\tau_{\rho z_m,l} \in K[X]  \qquad (\rho\in Re(m),\; 1\le m\le t,\; 1\le l\le d)$$
by \zitat{\hierarchies, Section 3}.  

Some of these polynomials are immediately recognizable. Namely, if $\rho\in Re(m)$ and $\rho z_m$ happens to be a $\S$-critical path, then $\rho z_m \seq \sum_{b\in\S(\rho z_m)} X_{\rho z_m,b}\,b$ and so
$$\tau_{\rho z_m,l} = \cases X_{\rho z_m,b_l}  &(\text{if\;} b_l\in \S(\rho z_m))\\  0  &(\text{if\;} b_l\notin \S(\rho z_m)) \endcases.  \tag4.1$$

Algorithmically, the $\seq$-expansion in terms of $\S$ of any path $pz_r \in \Phat$ of length $\le L$ is carried out by way of the following successive steps ($L$ steps at the most).  If $p z_r \in \S$, we are done.  Otherwise, suppose that $p z_r = p_2 p_1 z_r$, where $p_1 z_r$ is the unique $\S$-critical  initial subpath of $p z_r$.  Then  $p z_r\seq \sum_{b_l \in \S(p_1 z_r)} X_{p_1 z_r, b_l}\, p_2 b_l$.   If all of the paths $p_2 b_l$ belong to $\S$, we are done.  Otherwise, we observe that the unique $\S$-critical initial subpath of any path $p_2 b_l \notin \S$ is longer than $p_1 z_r$.  We repeat the procedure for the $p_2 b_l$ outside $\S$.  After a finite number of iterations, the process will thus lead to a $K[X]$-linear combination of paths in $\S$ and paths of length $> L$ in $\Phat$.  Since all of the latter paths are $\seq 0$, this will terminate the substitution process. 

As for decomposing $\grassS$ into its irreducible components:  There is a standard computer package for obtaining polynomials for the irreducible components of $\grassS$ from the $\tau_{\rho z_m,l}$;  it has been attached to the computer code for the calculation of the $\tau_{\rho z_m,l} \in K[X]$ in \zitat{\codes}. 

For the purpose of illustration, we will repeatedly refer to Carlson's Example, the local self-injective biserial algebra$$A= KQ/ \langle \, \alpha^2,\, \beta^2,\, \alpha\beta-\beta\alpha\, \rangle \,, \quad \text{where}\quad Q= \vcenter{\xymatrix{ 1 \ar@(ul,dl)_{\alpha} \ar@(ur,dr)^{\beta} }}.$$

\definition{Example 4.1} 
Let $A$ be the Carlson algebra and $\S$ the skeleton with graph
$$\xymatrixcolsep{1.25ex}\xymatrixrowsep{1.5pc}\xymatrix{
 &1 \dropup{z_1} \edge[dl]_{\alpha} \edge[dr]^{\beta} &&&1 \dropup{z_2} \edge[d]^{\alpha} &&1 \dropup{z_3} \drbl  \\
1 \edge[d]_{\beta} &&1  &&1  \\
1
}$$
in $\Phat$, where $P= \bigoplus_{1\le r\le 3} Az_r \cong (Ae_1)^3$. The $\S$-critical paths are
$$\alpha^2z_1\,,\; \alpha\beta z_1\,,\; \beta^2z_1\,,\; \beta z_2\,,\; \alpha^2 z_2\,,\; \beta\alpha z_2\,,\; \alpha z_3\,,\; \beta z_3\,,$$
and a set of generators for the left ideal $I\subseteq KQ$ is
$$R= \{ \alpha^2,\, \beta^2,\, \alpha\beta-\beta\alpha \} \cup \{ \text{the paths of length\;} 3 \}.$$
Since $L=2$ and all paths of length $\ge3$ are $\seq0$ in $(KQ[X]\otimes_{KQ} \Phat)\, / \, \U$, the latter paths do not yield any nonzero polynomials for $\grassS$, whence they will be ignored.

In our example,
$$\Rhat= \{ \alpha^2z_r\,,\, \beta^2z_r\,,\, \alpha\beta z_r- \beta\alpha z_r \mid 1\le r\le 3\},$$
and the basic equivalences in $KQ[X]\otimes_{KQ} \Phat$ modulo $\U$ are
$$\align
\alpha^2 z_1 &\seq X_{\alpha^2z_1,\beta\alpha z_1} \beta\alpha z_1\,, \qquad\qquad\qquad\qquad\qquad\qquad \alpha\beta z_1 \seq X_{\alpha\beta z_1,\beta\alpha z_1} \beta\alpha z_1  \\
\beta^2z_1 &\seq X_{\beta^2z_1,\beta\alpha z_1} \beta\alpha z_1  \\
\beta z_2 &\seq X_{\beta z_2,\alpha z_1}\alpha z_1+ X_{\beta z_2,\beta z_1} \beta z_1+ X_{\beta z_2,\beta\alpha z_1} \beta\alpha z_1+ X_{\beta z_2,\alpha z_2} \alpha z_2  \\
\alpha^2z_2 &\seq X_{\alpha^2z_2,\beta\alpha z_1} \beta\alpha z_1\,, \qquad\qquad\qquad\qquad\qquad\qquad \beta\alpha z_2 \seq X_{\beta\alpha z_2,\beta\alpha z_1} \beta\alpha z_1  \\
\alpha z_3 &\seq X_{\alpha z_3,\alpha z_1} \alpha z_1+ X_{\alpha z_3,\beta z_1} \beta z_1+ X_{\alpha z_3,\beta\alpha z_1} \beta\alpha z_1+ X_{\alpha z_3,\alpha z_2} \alpha z_2  \\
\beta z_3 &\seq X_{\beta z_3,\alpha z_1} \alpha z_1+ X_{\beta z_3,\beta z_1} \beta z_1+ X_{\beta z_3,\beta\alpha z_1} \beta\alpha z_1+ X_{\beta z_3,\alpha z_2} \alpha z_2 \,.
\endalign$$
Since $\alpha^2z_1$, $\beta^2z_1$, $\alpha^2z_2$ are $\S$-critical, we immediately obtain the equations
$$X_{\alpha^2z_1,\beta\alpha z_1}= X_{\beta^2z_1,\beta\alpha z_1}= X_{\alpha^2z_2,\beta\alpha z_1} =0  \tag4.2$$
from (4.1). Expanding the remaining elements of $\Rhat$ in  $(KQ[X]\otimes_{KQ} \Phat)/\U$ and recording the resulting equations for $\grassS$, we next obtain 
$$X_{\alpha\beta z_1,\beta\alpha z_1} =1  \tag4.3$$
from $\alpha\beta z_1- \beta\alpha z_1$ in $\Rhat$. The successive moves for the expansion of $\alpha\beta z_2- \beta\alpha z_2$ in $\Rhat$ are:
$$\align
\alpha\beta z_2- \beta\alpha z_2 &\seq X_{\beta z_2,\alpha z_1} \alpha^2z_1+ X_{\beta z_2,\beta z_1} \alpha\beta z_1+ X_{\beta z_2,\beta\alpha z_1} \alpha\beta\alpha z_1  \\
 & \qquad\qquad\qquad\qquad + X_{\beta z_2,\alpha z_2} \alpha^2z_2- X_{\beta\alpha z_2,\beta\alpha z_1} \beta\alpha z_1  \\
 &\seq X_{\beta z_2,\alpha z_1} X_{\alpha^2z_1,\beta\alpha z_1} \beta\alpha z_1+
 X_{\beta z_2,\beta z_1} X_{\alpha\beta z_1,\beta\alpha z_1} \beta\alpha z_1  \\
 & \qquad\qquad\qquad\qquad +
 X_{\beta z_2,\alpha z_2} X_{\alpha^2z_2,\beta\alpha z_1} \beta\alpha z_1-
 X_{\beta\alpha z_2,\beta\alpha z_1} \beta\alpha z_1 \,.
 \endalign$$
In light of (4.2) and (4.3), this results in the equation
 $$X_{\beta z_2,\beta z_1} - X_{\beta\alpha z_2,\beta\alpha z_1}=0.  \tag4.4$$
Expansion of the remaining four elements of $\Rhat$ gives us the equations
$$\aligned
X_{\beta z_2,\alpha z_1} &= - X_{\beta z_2,\alpha z_2} X_{\beta\alpha z_2,\beta\alpha z_1}  \\
X_{\alpha z_3,\alpha z_1} &= X_{\beta z_3,\beta z_1}- X_{\alpha z_3,\alpha z_2} X_{\beta\alpha z_2,\beta\alpha z_1}  \\
X_{\alpha z_3,\beta z_1} &= 0  \\
X_{\beta z_3,\alpha z_1} &= - X_{\beta z_3,\alpha z_2} X_{\beta\alpha z_2, \beta\alpha z_1}\,. 
\endaligned  \tag4.5$$
In summary, $\grassS \subseteq \AA^N$ is pinned down by the above equations (4.2--4.5) in the coordinates of $\AA^N$.
\enddefinition

A first theoretical application of the polynomials for $\grassS$ will be given in the second part of the upcoming section; see Proposition 5.2 and Theorem 5.4.

For a skeleton $\S$ in $\bPhat$, the process described above -- with $\bPhat$ taking the place of $\Phat$ and the index set $N$ enlarged to accommodate an increased number of $\S$-crit\-ical paths -- yields polynomials $\tau_{\rho z_m,l}$ in the correspondingly enlarged polynomial ring $K[X]$.  Their  simultaneous vanishing locus is $\biggrassS$ \zitat{\hierarchies, Section 3.D}. It turns out that the ``new" variables make no appearance in the $\tau_{\rho z_m,l}$, which yields the above-mentioned isomorphism $\biggrassS \cong \grass(g. \S) \times \AA^m$, where $g$ is a permutation automorphism in $\aut_A(\bPhat)$ taking  $\S$ to a subset of $\Phat$; see \zitat{\hierarchies, Theorem 3.17}.

\head 5. Carlson's Example revisited and generalized \endhead

For the first part of this section, we let $A$ be the Carlson algebra, as in Example 4.1.   (Carlson used this algebra to illustrate failure of cancellation for the degeneration order.  His argument is based on a comparison of orbit dimensions of objects in $\replad$;  see \zitat{\Rie, Section 3.1}.)

The irreducible components of the varieties $\replad$ are given in  \zitat{\RRS, Theorem 2.3} in terms of the generic structure of their modules. To make our technique for computing projective equations for the irreducible
components of the varieties $\GRASS(A,d)$ readily applicable to the Carlson
algebra, we list the generic radical layerings, together with their generic
skeleta. They are as follows:
\smallskip

{\bf $d$ even:} $\SS= (S_1^{a+b},S_1^{2a+b},S_1^a)$ such that $a,b\ge0$ with $4a+2b=d$, and a corresponding generic skeleton is
$$\xymatrixcolsep{1.25ex}\xymatrix{
 &1\dropup{z_1} \edge[dl]_{\alpha} \edge[dr]^{\beta} &&&&1\dropup{z_2} \edge[dl]_{\alpha} \edge[dr]^{\beta} &&&&&&1\dropup{z_a} \edge[dl]_{\alpha} \edge[dr]^{\beta} \edge[dr]^{\beta} &&&1\dropup{z_{a+1}} \edge[d]_{\alpha} &&&&1\dropup{z_{a+b}} \edge[d]_{\alpha} \\
 1 \edge[dr]_{\beta} &&1 &, &1 \edge[dr]_{\beta} &&1 &, &\cdots &, &1 \edge[dr]_{\beta} &&1 &, &1 &, &\cdots &, &1  \\
 &1 &&&&1 &&&&&&1
 }$$
 \smallskip
{\bf $d$ odd:} $\SS= (S_1^{a+b+1},S_1^{2a+b},S_1^a)$ for $a,b\ge0$, and $\SS= (S_1^{a+b},S_1^{2a+b+1},S_1^a)$ for $a\ge0$, $b\ge1$, where $4a+2b+1=d$. Corresponding generic skeleta are
$$\xymatrixcolsep{1.25ex}\xymatrix{
 &1\dropup{z_1} \edge[dl]_{\alpha} \edge[dr]^{\beta} &&&&1\dropup{z_2} \edge[dl]_{\alpha} \edge[dr]^{\beta} &&&&&&1\dropup{z_a} \edge[dl]_{\alpha} \edge[dr]^{\beta} \edge[dr]^{\beta} &&&1\dropup{z_{a+1}} \edge[d]_{\alpha} &&&&1\dropup{z_{a+b}} \edge[d]_{\alpha} &&1\dropup{z_{a+b+1}}  \\
  1 \edge[dr]_{\beta} &&1 &, &1 \edge[dr]_{\beta} &&1 &, &\cdots &, &1 \edge[dr]_{\beta} &&1 &, &1 &, &\cdots &, &1 &,  \\
 &1 &&&&1 &&&&&&1
 }$$
and
$$\xymatrixcolsep{1.25ex}\xymatrix{
 &1\dropup{z_1} \edge[dl]_{\alpha} \edge[dr]^{\beta} &&&&1\dropup{z_2} \edge[dl]_{\alpha} \edge[dr]^{\beta} &&&&&&1\dropup{z_a} \edge[dl]_{\alpha} \edge[dr]^{\beta} \edge[dr]^{\beta} &&&1\dropup{z_{a+1}} \edge[d]_{\alpha} &&&&1\dropup{z_{a+b-1}} \edge[d]_{\alpha} &&1\dropup{z_{a+b}} \edge[d]_{\alpha} \edge[dr]^{\beta}  \\
  1 \edge[dr]_{\beta} &&1 &, &1 \edge[dr]_{\beta} &&1 &, &\cdots &, &1 \edge[dr]_{\beta} &&1 &, &1 &, &\cdots &, &1 &, &1 &1  \\
 &1 &&&&1 &&&&&&1
 }$$
respectively.

\medskip

Via the technique of Section 4, it is easily verified that, over the Carlson algebra, all of the varieties $\grassS$ are full affine spaces.  In fact, this readily generalizes to arbitrary self-injective algebras with $J^3 = 0$.  Once justified, the upcoming theorem will bring the transfer results of Section 3 to bear.  

We return to the general notational conventions of Section 1. For the remainder of Section 5, we focus on self-injective algebras with vanishing radical cube.  Given any module $M$ over such an algebra $A$, we can clearly decompose $M$ in the form $M= M_0\oplus M_1$ where $J^2M_0 =0$ and $M_1$ is a direct sum of indecomposable projective modules of Loewy length $3$; indeed, by self-injectivity of $A$, the indecomposable $A$-modules of Loewy length $3$ are projective-injective. Accordingly, we call a skeleton $\S$ {\it normalized\/} if it is a disjoint union $\S= \S_0 \sqcup \S_1$, where $S_0$ consists of trees with edge paths of lengths at most $1$ and $\S_1$ is a union of skeleta of projectives with Loewy length $3$. The mentioned decompositions of $A$-modules makes the following lemma immediate.

\proclaim{Lemma 5.1} Suppose $A$ is self-injective with $J^3=0$. Then $\grassbd$ is covered by $\biggrassS$ where $\S$ traces the normalized skeleta with $\underbardim \S = \bd$. 
\qed\endproclaim

\proclaim{Proposition 5.2} Suppose $A$ is self-injective with $J^3=0$. Then the varieties $\grassS$, where $\S$ runs through the normalized skeleta over $A$, are full affine spaces.

Consequently, all $\grassbd$ are covered by varieties $\biggrassS$ isomorphic to full affine spaces.
\endproclaim

Instead of giving a formal argument resting on the analysis of a system of equations as in Section 4, we describe the setup leading to such a system and indicate how to proceed by induction.  For the remainder of Section 5, we let $A$ be a self-injective algebra with $J^3=0$, and $A_0= A/J^2A$ the factor algebra with vanishing radical square. According to Lemma 5.1, we may restrict our attention to normalized skeleta. Letting $\S= \S_0 \sqcup \S_1$ be a normalized skeleton, we observe that  each tree in $\S_1$ has the form
$$\xymatrixcolsep{2pc}\xymatrixrowsep{3pc}
\vcenter
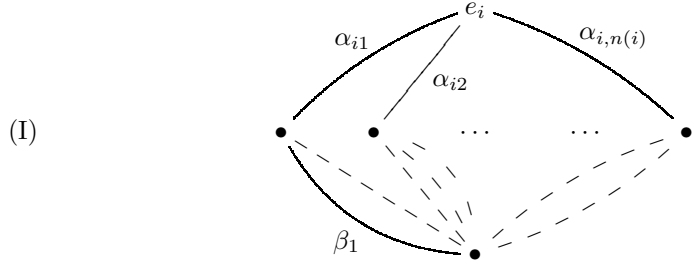
{\xymatrix{
 &&e_i \edge@/_/[dll]_{\alpha_{i1}} \edge[dl]^{\alpha_{i2}} \edge@/^/[drr]^{\alpha_{i,n(i)}}  \\
\bullet \edge@/_3ex/[drr]_{\beta_1} \dashedge[drr]   &\bullet \dashedge[dr] \dashedge@/^/[dr]  \dashedge@/^3ex/[dr]  &\cdots&\cdots &\bullet \dashedge@/_/[dll] \dashedge@/^/[dll]  \\
 &&\bullet
}} \tag \text{I}$$
\noindent Here $\alpha_{i1},\dots,\alpha_{i,n(i)}$ are all the arrows which start in the vertex $e_i$.  The possible values on $\grassS$ of the variables indexed by the corresponding $\S$-critical paths (the latter are indicated by dashed lines) are completely determined by the structure constants of $A$. 

As for the  trees in $\S_0$: They all have Loewy length $1$ or $2$, where the Loewy length equals $1$ plus the length of a longest edge path. It is convenient to induct on the number of trees of Loewy length $2$, these being of the form
$$\xymatrixcolsep{2pc}\xymatrixrowsep{3pc}
\vcenter
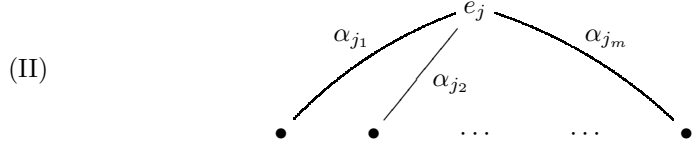
{\xymatrix{
 &&e_j \edge@/_/[dll]_{\alpha_{j_1}} \edge[dl]^{\alpha_{j_2}} \edge@/^/[drr]^{\alpha_{j_m}}  \\
\bullet & \bullet &\cdots&\cdots & \bullet
}} \tag \text{II}$$
where $\alpha_{j_1},\dots,\alpha_{j_m}$ is any list of (some) distinct arrows starting in $e_j$. 

The relations in $R$ which are relevant in the substitution process for obtaining $\grassS$ all have the form
$$p- c \, q \qquad \text{for distinct paths\ } p, q \in KQ  \text{\ and scalars\ } c \in K\, ,$$
with $c =0$ being permissible.  Analyzing the resulting system of equations by induction on the number of trees of type (II) in $\S$, one finds:  The variables $X_{b'_i, b_j}$ for $(b_i', b_j) \in N$ may be partitioned into 
\roster
\item"(a)" variables which are constant on $\grassS$;
\item"(b)" variables which, independently, attain arbitrary $K$-values on $\grassS$; and
\item"(c)" variables which arise in equations
$$X_{b'_i, b_j}= p_{ij\, , }$$
where $p_{ij}$ is a polynomial in variables belonging to groups (a) and (b).
\endroster

\proclaim{Lemma 5.3} Let $A$ be self-injective with $J^3=0$, and suppose $\S$ and $\tau$ are normalized skeleta with $\SS(\S)= \SS(\tau)$ such that $\biggrassS$ and $\GRASS(\tau)$ are both nonempty. Then $\biggrassS \cap \GRASS(\tau) \ne \varnothing$.
\endproclaim

\demo{Proof} Given that $\S$ and $\tau$ are normalized, we write  $\S= \S_0\sqcup \S_1$ and $\tau= \tau_0\sqcup \tau_1$ as above. The assertion amounts to the existence of an $A$-module $Z$ with the property that both $\S$ and $\tau$ are skeleta of $Z$. In view of the mentioned decomposition property of the objects in $\Amod$, there is precisely one $A$-module $Z_1$, up to isomorphism, which has skeleton $\S_1$, and our hypotheses guarantee that $Z_1$ also has skeleton $\tau_1$. In view of $\SS(\S)= \SS(\tau)$, we therefore obtain $\SS(\S_0)= \SS(\tau_0)$. Nontriviality of $\biggrassS$ and $\GRASS(\tau)$ moreover guarantees nontriviality of $\GRASS(\S_0)$ and $\GRASS(\tau_0)$. Since $\S_0$ and $\tau_0$ are skeleta for modules over the truncated path algebra $A_0$, the intersection $\GRASS(\S_0) \cap \GRASS(\tau_0)$ is nonempty by \zitat{\BHT, Theorem 5.3}. In other words, $\S_0$ and $\tau_0$ are both skeleta of some $A_0$-module $Z_0$. Thus, the direct sum $Z := Z_0\oplus Z_1$ is as required.
\qed\enddemo

Lemma 5.3 implies that all of the varieties $\biggrassSS$ are irreducible. Therefore, the results of Section 3 now yield the following consequences:

\proclaim{Theorem 5.4}  Again, suppose that $A$ is self-injective with $J^3 = 0$, and let $\SS$ be any semisimple sequence.

The varieties $\biggrassSS$ and $\repSS$ are irreducible, unirational, and smooth, and the irreducible components of $\grassbd$ {\rm{(}}resp., $\replabd${\rm{) }}are among the closures of the $\biggrassSS$ {\rm{(}}resp., $\repSS${\rm{)}} in $\grassbd$ {\rm{(}}resp., $\replabd${\rm{)}}. 

In particular,  all irreducible components of $\grassbd$ and $\replabd$ are unirational.
Moreover, given any such irreducible component with generic radical layering $\SS$, all of its singularities lie in 
$$\overline{\biggrassSS} \setminus \biggrassSS\qquad \text{resp.} \qquad \overline{\repSS} \setminus \repSS.$$
\endproclaim

\head 6. Projective equations \endhead

In situations where the irreducible components of $\grassbd$ are known --  there are few so far  --   the generic skeleta of the modules parametrized by the individual components are immediate from the  more detailed representation-theoretic information inherent in the classification.  This is the piece of information -- still missing aside from a few cases -- which makes the following general method for computing polynomials for the closures of the components of  $\biggrassSS$ in the ambient projective space $\PP$ applicable towards the more specific task of finding polynomial equations for the irreducible components of $\grassbd$ in $\PP$.   

We know that the irreducible components of $\biggrassSS$ are in natural correspondence with the irreducible components of the $\biggrassS$, where $\S$ traces the skeleta with radical layering $\SS$ (see Section 2.C).
Thus, we first indicate how polynomials for $\grassS$ determine polynomials for $\biggrassS$; this is simply a matter of bookkeeping. 
Subsequently, in three steps, we describe how to algorithmically compute equations which determine the closures of irreducible components of $\biggrassS$ in the projective space $\PP$. 

In the following, we assume that $\grassS$ is nonempty, which is equivalent to $\biggrassS \ne \varnothing$.

\definition{From $\grassS$ to $\biggrassS$}
Recall that $d= |\bd|= \sum_{i=1}^n d_i$. Without loss of generality, the skeleton $\S$ is contained in the distinguished projective cover $P= \bigoplus_{1\le r\le t} Az_r$ of the top $T$ of $\S$, with $\grassS$ located in $\PP(\la^{\dim P-d} P)$. The enlarged variety $\biggrassS$ lives in $\PP= \PP(\la^{\dim \bP-d} \bP)$, where $\bP= \bigoplus_{1\le s\le d} Az_s \supseteq P$; here $z_1,\dots,z_t$ and $z_1,\dots, z_d$ are the distinguished sequences of top elements of $P$ and $\bP$, respectively. As before, $\Phat$ and $\bPhat$ are the projective $KQ$-modules with the corresponding sequences of top elements.

In the enlarged setting of $\bPhat$, the set of $\S$-critical paths includes the paths $z_{t+1},\dots,z_d$ of length zero. On the other hand, all of the $\S$-critical paths in $\Phat$ have positive length. In \zitat{\hierarchies}, it was shown that any point $C\in \biggrassS$ can be uniquely written in the form
$$C = \sum \Sb b' \, \S \text{-critical}\\ \len(b')>0 \endSb A \biggl (b' \, - \, \sum_{b \in \S(b')} \cbb\, b \biggr) + \sum \Sb b' \, \S \text{-critical}\\ \len(b')=0 \endSb A \biggl (b' \, - \, \sum_{b \in \S(b')} \cbb\, b \biggr), \tag6.1$$
with no restrictions governing the scalars appearing in the right-hand sum. In other words, the polynomials $\tau_{\rho z_m,l}$ obtained in Section 4, now viewed as elements of a polynomial ring with additional variables $X_{b',b}$, where $b'$ is $\S$-critical of length $0$, determine $\biggrassS$ as well.   
  \enddefinition
 
\definition{Upgraded notational setup}  Write $a := \dim\bP-d$. 

In the sequel, we will systematically identify the skeleton $\S \subseteq \bPhat$ with a subset of $\bP$. Since $\biggrassS \ne \varnothing$, we know that $\S$ is a linearly independent subset of $\bP$ (indeed, $|\S|=d$, and $\S$ induces a basis for the $d$-dimensional module $\bP/C$, for any point $C\in \biggrassS$). It will be convenient to list $\S$ in the form 
$$\S= \{b_1,\dots,b_d\}.$$

Next, let $\sigma' \subseteq \bPhat$ be a set of $\S$-critical paths whose images in $\bP$ induce a basis for 
$$\biggl(\, \sum_{b'\, \S\text{-critical}} Kb' \,+\, \sum_{1\le j\le d} Kb_j \biggr) \bigg/ \biggl(\, \sum_{1\le j\le d} Kb_j \biggr) ;$$
say $\sigma'= \{b'_1,\dots,b'_u\}$.
Finally, supplement the linearly independent set $\S\sqcup \sigma'$ by suitable paths $b''_i$ in $\bP$ so that
$$B :=  \{b_1,\dots,b_d\} \sqcup \{b'_1,\dots,b'_u\} \sqcup \{b''_1,\dots,b''_v\}= \{b_1,\dots,b_d,b'_1,\dots,b'_u,b''_1,\dots,b''_v \}$$
is an ordered basis for $\bP$; note that $u+v=a$. We also write $B= \{w_1,\dots,w_{\dim\bP}\}$ in this ordering. 
By the choice of the $b_i'$, any $\S$-critical path $b'\notin\S'$ can be uniquely expressed as a $K$-linear combination
$$b'= \sum_{1\le j\le d} k_{\,b',j}\,b_j + \sum_{1\le i\le u} k'_{\,b',i}\,b'_i  \tag6.2$$
in $\bP$. Consequently, the expansion given in (6.1) reduces to
$$C = \sum_{b'\; \S\text{-critical}} A \biggl (b' \, - \, \sum_{b \in \S(b')} \cbb\, b \biggr)= \sum_{1\le i\le u} A \biggl( b'_i - \sum \Sb 1\le j\le d\\ b_j\in \S(b'_i) \endSb c_{ij} b_j \biggr),  \tag6.3$$
where $c_{ij} := c_{\,b'_i,b_j}$. 
This allows us to ignore the $b' \notin \S'$ and to move to the more convenient smaller index set
$$N_0 := \{ (i,j)\in \{1,\dots,u\}\times \{1,\dots,d\} \mid b_j\in \S(b'_i) \}.$$
Accordingly, write $X_{ij} := X_{b'_i,b_j}$ for the variables of the (potentially shrunken) polynomial ring $K[X]$.

The basis $\BB$ for the exterior power $\la^a\bP$ corresponding to the basis $B$ for $\bP$ is
$$\BB := \{ w_{i_1}\wedge\cdots\wedge w_{i_a} \mid (i_1,\dots,i_a) \in N_1 \},$$
where $N_1 := \{ (i_1,\dots,i_a) \mid 1\le i_1< \cdots< i_a\le \dim\bP \}$. We are thus viewing $\PP(\la^a\bP)$ as the $(|N_1|-1)$-dimensional projective space corresponding to the affine space $\AA^{N_1}$.
\enddefinition

\definition{Step 1. Affine equations for $\biggrassS$ and its components, viewed as closed subvarieties of $\AA^{N_0}$}
Recall the isomorphism from $\biggrassS$ onto a closed subvariety of $\AA^{\bold{N}}$ introduced at the end of Section 2.D.  Followed by the projection map onto $\AA^{N_0}$ along the coordinates in $\bold{N} \setminus N_0$, it yields an isomorphism from $\biggrassS$ onto a closed subset of $\AA^{N_0}$.   This is due to the fact that, in light of (6.2) and (6.3), we may discard the variables $X_{b',b_j}$ for $\S$-critical paths $b'\notin \S'$ without any loss of information.  Indeed, one checks that
$$X_{b',b_j}= \cases k_{\,b',j}+ \sum_{1\le i\le u} k'_{\,b',i}X_{ij}  &(b_j\in \S(b'))\\  k_{\,b',j} &(b_j\notin \S(b')) \endcases \qquad\quad (1\le j\le d)  \tag6.4$$
for such paths $b'$, where the scalars $k_{\sssize\bullet, \bullet}$ and $k'_{\sssize\bullet, \bullet}$ are as in (6.2). 
The image of the isomorphism which assigns to any $C \in \biggrassS$ the point $(c_{ij}) \in \AA^{N_0}$ which is determined by the requirement (6.3) will again be denoted by $\biggrassS$; it differs from the incarnation in $\AA^{\bold{N}}$ only by the disappearance of redundant coordinates.  One obtains it as the zero locus of the polynomials
$$\tau_{\rho z_m,l} \in K[X]= K[X_{ij} \mid (i,j) \in N_0 ] \qquad\quad (\rho\in R\,e(m),\; 1\le m\le t,\; 1\le l\le d)$$
as described in Section 4, with substitutions (6.4) as needed.
\enddefinition

\definition{An illustration of Step 1} We will use the Carlson algebra $A$ discussed in Section 5 as a running example to illustrate the procedures being presented.  For $d = 2$, the variety $\GRASS(A,d)$ is clearly isomorphic to $\GRASS(A_0,d)$, where 
$$A_0= KQ/ \langle \, \alpha^2,\, \beta^2,\, \alpha\beta,\, \beta\alpha\, \rangle \quad \text{with}\quad Q= \vcenter{\xymatrix{ 1 \ar@(ul,dl)_{\alpha} \ar@(ur,dr)^{\beta} }}.$$
  The only generic radical layering for $\GRASS(A_0, 2)$ is $\SS = (S_1, S_1)$, and therefore $\GRASS(A_0,2)$ equals the closure $\overline{\biggrassSS}=\overline{\biggrassS}$ for any skeleton $\S$ with radical layering $\SS$.
Let us consider the skeleton $\S= \{z_1,\alpha z_1\} \subset \bP= A_0z_1 \oplus A_0z_2$, together with $b'_1= \beta z_1$ and $b'_2= z_2$, following the notation above.  From Proposition 5.2 we know that $\grassS$ and $\biggrassS$ are full affine spaces in this case. In fact, with respect to the coordinates $X_{ij}$ for $\AA^{N_0}$, the algorithm of Section 4 does not yield any equations, that is, $\biggrassS= \AA^{N_0}$.
\enddefinition

\definition{Back to the general case} In the sequel, we focus on a fixed irreducible component $\C(\S)$ of $\biggrassS$ and compute polynomials for its closure in $\PP= \PP(\la^a\bP)$ -- equivalently, for the closure of $\C(\S)$ in $\grassbd$  -- in two further steps. The variables $Y_b$ ($b\in\BB$) for this projective space will be abbreviated in the following form:
$$Y_{i_1,\dots,i_a} = Y_{w_{i_1}\wedge \cdots\wedge w_{i_a}} \,,$$
where $(i_1,\dots,i_a) \in N_1$ (recall the notation set up ahead of Step 1). 
In particular, $Y_{d+1,\dots,d+a}= Y_{b'_1\wedge \cdots\wedge b'_u\wedge b''_1\wedge \cdots\wedge b''_v}$. Thus, the homogeneous coordinate ring of $\PP$ is the ring
$$K[Y] := K[Y_b\mid b\in \BB]= K[Y_{i_1,\dots,i_a} \mid (i_1,\dots,i_a) \in N_1 ].$$
It will, moreover, be convenient to abbreviate the variable $Y_{i_1,\dots,i_a}$ representing
$$w_{i_1}\wedge \cdots\wedge w_{i_a} = b_l\wedge b'_1\wedge \cdots b'_{k-1}\wedge b'_{k+1}\wedge \cdots\wedge b'_u \wedge b''_1\wedge \cdots\wedge b''_v$$
 for $(k,l)\in N_0$ by $\hatY_{kl}$.

We will find that the homogenizations (relative to the variable $Y_{d+1,\dots,d+a}$) of the polynomials for $\C(\S)$ are among the polynomials for the closure, and that the additional homogeneous polynomials in $K[Y]$ will arise as ``trivial inflations" of the polynomials we already have.
\enddefinition

\definition{Step 2.  Affine equations for the relative closure of $\C(\S)$ in the Schubert cell $\SchubertS$ of the Grassmannian $\Gr(a,\bP) \subseteq \PP$} Recall from Section 2.B that $\SchubertS$  is the dehomogenization of  $\Gr(a,\bP)$ relative to the coordinate $Y_{b'_1\wedge \cdots\wedge b'_u\wedge b''_1\wedge \cdots\wedge b''_v} = Y_{d+1,\dots,d+a}$.  We note  that $\C(\S) \subseteq \SchubertS$, and let $\overline{\C(\S)}$ be the closure of $\C(\S)$ in $\PP$.

On the road towards developing affine polynomials for the irreducible variety $\overline{\C(\S)}\cap \SchubertS$ in the big Schubert cell $\SchubertS$ of the Grassmannian $\Gr(a,\bP)$, we first observe that this intersection is nothing new.  
\enddefinition

\proclaim{Lemma 6.1}  $\overline{\C(\S)}\cap \SchubertS = \C(\S)$.
\endproclaim

\demo{Proof}  Let $\SS = (\SS_1, \dots, \SS_L)$ be the radical layering of $\S$, that is,  $\SS_l = \bigoplus_{1 \le i \le n} S_i^{m_{li}}$, where $m_{li}$ is the number of those paths of length $l$ in $\S$ which end in the vertex $e_i$.  

{\bf Claim:} $\SchubertS \cap \grassbd \subseteq  \bigcup_{\SS' \le \SS} \biggrass(\SS')$, where $\le$ is the dominance order of \zitat{\hierarchies, Definition 2.10}.

Indeed, let $C$ be in the above intersection. Then $\{b_i+C \mid 1\le i\le d\}$ is a basis for $\bP/C$. Since $b_i+C \in J^l(\bP/C)$ whenever $b_i$ has length $l$, we thus obtain
$$\bigoplus_{0\le l\le m} \SS(\bP/C) \le \bigoplus_{0\le l\le m} \SS_l$$
for $0\le m\le L$, that is, $\SS(\bP/C) \le \SS$. This proves the claim.

On the other hand, $\overline{\biggrassSS} \subseteq \bigcup_{\SS' \ge \SS} \biggrass(\SS')$ \zitat{\hierarchies, Observation 2.11}, and because 
$\C(\S) \subseteq \biggrassS \subseteq \biggrassSS$, we conclude that $\overline{\C(\S)} \cap \SchubertS \subseteq \biggrassSS$. Since $\SchubertS \cap \biggrassSS = \biggrassS$, this completes the argument.
\qed \enddemo

Hence, our task consists of finding polynomial equations for $\C(\S)$ in the alternate coordinate system of $\SchubertS$ induced by the basis $\BB$ for $\la^a(\bP)$ (cf\. the notation ahead of Step 1). We prepare with a preliminary shift from the incarnation of $\biggrassS$ in the affine space $\AA^{N_0}$ to $\la^a\bP$, prior to expressing the resulting points in $\la^a\bP$ in the coordinates for $\SchubertS$ determined by $\BB$. Namely, we first specify an assignment
$$C= (c_{ij}) \in \AA^{N_0} \quad\longmapsto\quad \bc= C'_1\wedge \cdots\wedge C'_u \wedge C_1'' \wedge \cdots\wedge C_v'' \in \la^a\bP.$$

\proclaim{Lemma 6.2} There exist polynomials
$q_{kl} \in K[X]= K[X_{ij} \mid (i,j)\in N_0]$
for $1\le k\le v$ and $1\le l\le d$ {\rm(}depending only on our choice of ordered basis $B$ for $\bP${\rm)} with the following property:

Any point $C\in \biggrassS$ -- say $C$ is given by $(c_{ij})_{(i,j)\in N_0}$ in $\AA^{N_0}$ -- is the projective point in $\SchubertS \subseteq \PP(\la^a\bP)$ which is represented by the element
$$\bc= C'_1\wedge \cdots\wedge C'_u\wedge C''_1\wedge \cdots\wedge C''_v$$
in the exterior power $\la^a\bP$, where
\roster
\item $C'_i= b'_i- \sum_{1\le j\le d,\; (i,j)\in N_0} c_{ij}b_j \in \bP$ for $1\le i\le u$;
\item $C''_k = b''_k- \sum_{l=1}^d q_{kl}(c_{ij})b_l \in \bP$ for $1\le k\le v$.
\endroster
In particular, the $Y_{d+1,\dots,d+a}$-coordinate {\rm(}$=$ the $Y_{b'_1\wedge \cdots\wedge b'_u\wedge b''_1\wedge \cdots\wedge b''_v}$-coordinate{\rm)}  of this point $\bc$ is $1$.
\endproclaim

\demo{Proof} As $C= \sum_{i=1}^u A\,C'_i$, the congruences $b'_i \equiv \sum_j c_{ij}b_j \pmod C$ permit us to expand the $b''_k$, modulo $C$, in terms of $b_1,\dots,b_l$ for $1\le k\le v$, which results in elements $C_k''$ of type (2) in $C$. (Each $b''_k$ has the form $pb'$, where $b'$ is the unique $\S$-critical initial subpath of $b''_k$ and $p$ is a path of positive length, and we proceed by induction on $\len(p)$.) In this expansion, the coefficients of the $b_l$ are polynomial expressions in the $c_{ij}$. These expansions give rise to the $q_{kl} \in K[X]$ as postulated.
\qed\enddemo

The following lemma is an immediate consequence of Lemma 6.2. It supplements the shift $C= (c_{ij}) \in \AA^{N_0} \mapsto \bc\in \la^a\bP$ by specifying the affine coordinates of $\bc$ in the standard coordinatization of the Schubert cell $\SchubertS$ of $\Gr(a,\bP)$ relative to $\BB$. Namely, given that this cell consists of the points in $\Gr(a,\bP)$ whose $Y_{d+1,\dots,d+a}$-coordinates are nonzero, $\SchubertS$ is identified with the intersection of $\Gr(a,\bP)$ and the hyperplane $Y_{d+1,\dots,d+a}=1$ in $\AA^{N_1}$.  Any point $\bc$ as above belongs to this hyperplane by Lemma 6.2.  Moreover, it is immediate that the coefficient of $C$ corresponding to the basis element $\hatY_{ij}$ equals $\pm c_{ij}$ for $(i,j)\in N_0$. We pin down the sign change by $\vareps_{ij}\in \{1,-1\}$: Thus, the $\hatY_{ij}$-coordinate of $\bc$ in $\la^a\bP$ is $\vareps_{ij}c_{ij}$. 

\proclaim{Lemma 6.3} There exist polynomials
$$\rho_{i_1,\dots,i_a} \in K[X]  \qquad \text{for\ } (i_1,\dots,i_a) \in N_1 \,,$$
depending only on our choice of ordered basis $B$ for $\bP$,  such that the following hold:

$\bullet$ Any point $C\in \biggrassS$, with coordinates $(c_{ij})_{(i,j)\in N_0}$ in $\AA^{N_0}$, is represented by the element
$$\sum_{(i_1,\dots,i_a) \in N_1} \rho_{i_1,\dots,i_a}(c_{ij})\,\ w_{i_1}\wedge \cdots\wedge w_{i_a}$$
relative to the basis $\BB$ for $\la^a\bP$, and the expansion coefficient $\rho_{d+1,\dots,d+a}(c_{ij})$ of the basis element $b'_1\wedge \cdots\wedge b'_u\wedge b''_1\wedge \cdots\wedge b''_v$ is the constant $1$.

$\bullet$ If, for $(k,l)\in N_0$, we write $\rhohat_{kl}$ for the expansion polynomial
$$\rho_{l,d+1,\dots,d+k-1,d+k+1,\dots,\dim\bP}$$
corresponding to $\;b_l\wedge b'_1\wedge \cdots b'_{k-1}\wedge b'_{k+1}\wedge \cdots\wedge b'_u \wedge b''_1\wedge \cdots\wedge b''_v\;$, then
$$\rhohat_{kl}= \vareps_{kl}X_{kl}$$
where $\vareps_{kl}\in \{\pm1\}$ is as introduced above.
\qed\endproclaim

\proclaim{Theorem 6.4} The affine variety $\overline{\C(\S)}\cap \SchubertS = \C(\S)$ is determined by the following polynomials in the variables $Y_{i_1,\dots,i_a}$ for $\SchubertS$: 
$$Y_{d+1,\dots,d+a}-1 \,, \qquad Y_{i_1,\dots,i_a} - \rho_{i_1,\dots,i_a}(\vareps_{ij} \hatY_{ij}) \,, \qquad \tau_{\rho z_m,l}(\vareps_{ij} \hatY_{ij}) \,, \tag6.5$$
where the $ \rho_{i_1,\dots,i_a} \in K[X]$ are as specified in Lemma {\rm 6.3}, and the $\tau_{\rho z_m,l} \in K[X]$ are the polynomials for $\biggrassS$ determined in Step {\rm 1}.
\endproclaim

\demo{Proof}  By Lemma 6.1, $\C(\S)$ is a closed subvariety of $\AA^{N_1}$. Let $V \subseteq \AA^{N_1}$ be the zero locus of the polynomials (6.5) in $K[Y]$.
From Lemma 6.3, we derive that $\C(\S) \subseteq V$.

For the reverse inclusion, it suffices to observe that, given any point $D\in V$, the point
$$C:= (\vareps_{ij}D_{ij})_{(i,j)\in N_0} \in \AA^{N_0}$$
belongs to the incarnation of $\C(\S)$ in the affine coordinatization of Step 1. Here $D_{ij}$ is the $\hatY_{ij}$-coordinate of $D$. It follows from our construction that the remaining $\AA^{N_1}$-coordinates of this point $C$ coincide with those of $D$, and therefore $D\in \C(\S) \subseteq \AA^{N_1}$.
\qed\enddemo

The proof of Theorem 6.4 exhibits, in particular, a ``sign-weighted projection'' 
$$\SchubertS \supseteq \C(\S) \ni D \quad\longmapsto\quad (\vareps_{ij}D_{ij})_{(i,j)\in N_0} \in \AA^{N_0},$$
the image of which is the incarnation of $\C(\S)$ in $\AA^{N_0}$.

\definition{The example: Step 2 for $\GRASS(A_0,2)$} Choose the following ordered basis for $\bP$:
$$B := (b_1,b_2,b'_1,b'_2,b''_1,b''_2) = (z_1,\alpha z_1, \beta z_1, z_2, \alpha z_2, \beta z_2).$$
We compute the elements $C'_i$ and $C''_k$ appearing in Lemma 6.2 to be 
$$\xalignat2
C'_1 &= b'_1-c_{12}b_2  &C'_2 &= b'_2- c_{21}b_1- c_{22}b_2  \\
C''_1 &= b''_1- c_{21}b_2  &C''_2 &= b''_2- c_{12}c_{21}b_2 \,.
\endxalignat$$
Then the element $\bc \in \la^a\bP$ is
$$\multline
\bc = w_{3456} -c_{12}c^2_{21}w_{1235} +c^2_{21}w_{1236} -c_{12}c_{21}w_{1256}  +c_{21}w_{1356}  \\
+c_{12}c_{21}w_{2345} -c_{21}w_{2346}
  +c_{22}w_{2356} -c_{12}w_{2456} \,,
 \endmultline$$
where we communicate the basis vectors $w_{i_1}\wedge w_{i_2}\wedge w_{i_3}\wedge w_{i_4}$ for $\PP(\la^4\bP)$ in the form $w_{i_1i_2i_3i_4}$. Analogously abbreviating the polynomials $\rho_{i_1,\dots,i_4}$ of Lemma 6.3, we find that
$$\matrix
\rho_{3456} = 1  &\quad&\rho_{1235} = -X_{12}X^2_{21}  &\quad&\rho_{1236} = X^2_{21}  \\
\rho_{1256} = -X_{12}X_{21}  &&\rhohat_{21} = \rho_{1356}= X_{21}  &&\rho_{2345} = X_{12}X_{21}  \\
\rho_{2346} = -X_{21}  &&\rhohat_{22} = \rho_{2356}= X_{22}  &&\rhohat_{12} = \rho_{2456}= -X_{12} \,.  
\endmatrix$$
Consequently, Theorem 6.4 tells us that $\biggrassS \subseteq \SchubertS$ is the vanishing set of the following collection of polynomials:
$$\matrix
Y_{3456}-1  &\quad&Y_{1235}-\hatY_{12} \hatY_{21}^2  &\quad&Y_{1236}-\hatY_{21}^2  \\ 
Y_{1256}-\hatY_{12} \hatY_{21}  &&Y_{2345}+ \hatY_{12} \hatY_{21} &&Y_{2346}+ \hatY_{21}   \\
Y_{1234}  &&Y_{1245}  &&Y_{1246}  \\
Y_{1345}  &&Y_{1346}  &&Y_{1456} \,,
\endmatrix  \tag6.6$$
where $\hatY_{12}= Y_{2456}$, etc., according to the convention preceding Step 2.
\enddefinition

\definition{Step 3. Homogeneous polynomials for the closure $\overline{\C(\S)}$ in the projective space $\PP(\la^a\bP)$}
\enddefinition

\proclaim{Final Observation 6.5} {\rm(}We exclude the trivial case where $\bP$ is semisimple.{\rm)} Let
$\J \subseteq K[Y]$
be the ideal generated by the polynomials {\rm(6.5)} described in Theorem {\rm 6.4} of Step {\rm2}, and for any $\psi \in \J$, let $\psi^*$ be the homogenization of $\psi$ relative to the variable $Y_{d+1,\dots,d+a}$. Then $\overline{\C(\S)}$ is the zero locus of the homogeneous ideal $\J^*$ generated by the $\psi^*$ with $\psi\in \J$.
\endproclaim

\demo{Proof} This is routine.
\qed\enddemo

Comment: While the polynomials for $\C(\S)= \overline{\C(\S)}\cap \SchubertS$ obtained in Step 2 simply result from those computed in Step 1 by way of more expansive bookkeeping, the final homogenization step is best carried out by a computer. This calculation is rendered labor-intensive by the fact that the ideal generated by the polynomials
$$\bigl( Y_{i_1,\dots,i_a}- \rho_{i_1,\dots,i_a}(\vareps_{ij} \hatY_{ij}) \bigr)^* \qquad\text{and}\qquad \tau_{\rho z_m,l}(\vareps_{ij} \hatY_{ij})^* $$
alone may be too small to cut out $\overline{\C(\S)}$.

\definition{The example: Step 3 for $\GRASS(A_0,2)$} In the example, the homogenizations, with respect to $Z= Y_{3456}$, of the 12 polynomials displayed in (6.6) are not sufficient to determine the closure of $\biggrassS$ in $\PP(\la^4\bP)$. We claim that a suitable set of homogeneous polynomials is the following:
$$\matrix
Y_{1235}Z-\hatY_{12}Y_{1236}  &\quad&Y_{1236}Z-\hatY_{21}^2  &\quad&Y_{1256}Z-\hatY_{12} \hatY_{21}  \\
\hatY_{11}  &&Y_{2345}+Y_{1256}  &&Y_{2346}+ \hatY_{21}  \\
Y_{1234}  &&Y_{1245}  &&Y_{1246}  \\
Y_{1345}  &&Y_{1346} && \hatY_{12}Y_{1236}- \hatY_{21}Y_{1256}  \\
Y_{1236}Y_{1256}- \hatY_{21}Y_{1235}  &&Y^2_{1256}- \hatY_{12}Y_{1235} \,.
\endmatrix \tag6.7$$ 
These polynomials obviously vanish on $\biggrassS$, and the ideal $J$ they generate contains the homogenizations of the polynomials in (6.6). One can show (by hand in this small example, or with Macaulay2, for instance) that $J$ is a prime ideal of $K[Y]$, and therefore that $J$ is the vanishing ideal of $\overline{\biggrassS}$. In particular, the zero locus of the polynomials displayed in (6.7) is $\overline{\biggrassS}= \GRASS(A_0,2)$.

There is a unique point $C_0\in \GRASS(A_0,2)$ representing the semisimple module $S_1^2$, namely, the projective point in $\PP(\la^4\bP)$ corresponding to $b_2\wedge b'_1\wedge b''_1\wedge b''_2$. In terms of our given coordinates for $\PP(\la^4\bP)$, the point $C_0$ is given by $Y_{2356}=1$ and $Y_{i_1i_2i_3i_4} =0$ for all other $(i_1,\dots,i_4)$. Since $Y_{2356}= \hatY_{22}$ does not appear in (6.7), we see that the ring $K[Y]/( J+ \langle Y_{2356}-1\rangle)$ is not regular at the maximal ideal corresponding to $C_0$, meaning that $C_0$ is a singular point of $\GRASS(A_0,2)$.

Finally, we deduce normality of $\GRASS(A_0,2)$ from the above information. By \zitat{\Har, p\. 23, Exercise 3.18}, it suffices to show that the variety $\GRASS(A_0,2)$ is projectively normal, i.e., that its homogeneous coordinate ring $K[Y]/J$ is integrally closed. Observe that $K[Y]/J$ is a polynomial ring over $R_0/J_0$,  where
$$\align
R_0 &:= K[Z,\hatY_{12},\hatY_{21},Z_1,Z_2,Z_3]  \\
Z_1&:= Y_{1236},\, Z_2 := Y_{1256},\, Z_3 := Y_{1235}
\endalign$$
and $J_0$ is the ideal of $R_0$ generated by
$$\matrix
Z_3Z-\hatY_{12}Z_1  &&Z_1Z-\hatY_{21}^2  &&Z_2Z-\hatY_{12}\hatY_{21}  \\
\hatY_{12}Z_1-\hatY_{21}Z_2  &&Z_1Z_2-\hatY_{21}Z_3  &&Z_2^2-\hatY_{12}Z_3 \,.
\endmatrix$$
It is thus enough to show that $R_0/J_0$ is integrally closed (e.g., \zitat{\Mats, Proposition 17.B}).

One checks that the inclusion map $K[Z,\hatY_{12},\hatY_{21}] \rightarrow K[Z^{\pm1},\hatY_{12},\hatY_{21}]$ induces an embedding $\phi: R_0/J_0 \rightarrow K[Z^{\pm1},\hatY_{12},\hatY_{21}]$ such that
$$\xalignat3
\phi(Z_1+J_0) &= Z^{-1}\hatY_{21}^2  &\phi(Z_2+J_0) &= Z^{-1}\hatY_{12}\hatY_{21}  &\phi(Z_3+J_0) &= Z^{-2}\hatY_{12}\hatY_{21}^2 \,.
\endxalignat$$
(This is one route to show that $J_0$ is a prime ideal.)
It follows that $R_0/J_0$ is isomorphic to the semigroup algebra of the subsemigroup $Q\subseteq \ZZ^3$ generated by
$$(1,0,0),\, (0,1,0),\, (0,0,1),\, (-1,0,2),\, (-1,1,1),\, (-2,1,2).$$
It is not hard to see  that $Q$ is saturated in $\ZZ^3$, meaning that the intersection of $\ZZ^3$ with the cone $\RR_{\ge0}Q$ equals $Q$. Consequently, \zitat{\MS, Proposition 7.25} implies that $R_0/J_0$ is integrally closed, as desired.

We summarize. Recalling that $\GRASS(A,2) \cong \GRASS(A_0,2)$ and taking Theorems 3.3 and 5.4 into account, we conclude:
\roster
\item"" {\it For the Carlson algebra $A$, the varieties $\Rep(A,2)$ and $\GRASS(A,2)$ are irreducible, unirational and  normal, but not smooth.}
\endroster
\enddefinition
For $d>2$, the varieties $\replad$ and $\GRASS(A,d)$ are no longer irreducible, but their irreducible components are unirational by Theorem 5.4, and we expect these components to be normal. 

As mentioned in the introduction, we conjecture that, more generally, normality holds for all irreducible components of the varieties $\grassbd$ and $\replabd$, provided that $A$ is either a truncated path algebra or self-injective with vanishing radical cube.



\Refs
\widestnumber\key{{\bf 99}}

\ref\no \BHT \by E. Babson, B. Huisgen-Zimmermann, and R. Thomas\paper Generic representation theory of quiver with relations \jour J. Algebra \vol 322 \yr 2009 \pages 1877--1918 \endref

\ref\no \codes \bysame  \paper Maple codes for computing $\operatorname{\frak{Grass}}(\sigma)$'s \finalinfo posted at
www.math.washington.edu/\linebreak{}$\sim$thomas/programs/programs.html \endref

\ref \no \GeomIV  \by K. Bongartz and B. Huisgen-Zimmermann \paper Varieties of uniserial representations  IV. Kinship to geometric quotients \jour Trans. Amer. Math. Soc. \vol 353 \yr 2001 \pages 2091--2113  \endref

\ref\no \CBS \by W. Crawley-Boevey and J. Schr\"oer \paper Irreducible
components of varieties of modules \jour J. reine angew. Math. \vol 553 \yr
2002 \pages 201--220  \endref

\ref\no\djong \by A.J. de Jong, et al. \paperinfo The Stacks Project, www.math.columbia.edu/algebraic\underbar{\hphantom{x}}geometry/\linebreak{}stacks-git/ \endref

\ref\no \GP \by I. M. Gelfand and V. A. Ponomarev \paper Indecomposable representations of the Lorentz group \jour Russian Math. Surveys \vol 23 \yr 1968 \pages 1--58 \endref

\ref\no \EGAfourtwo \by A. Grothendieck and J. Dieudonn\'e \paper \'El\'ements de g\'eometrie alg\'ebrique IV. \'Etude locale des sch\'emas et des morphismes de sch\'emas, Seconde partie \jour Publ. math. I.H.E.S. \vol 24 \yr 1965 \pages 5--231 \endref

\ref\no \Har \by R. Hartshorne \book Algebraic Geometry \publaddr New York \yr 1977 \publ Springer-Verlag \endref

\ref\no \Hess \by W. Hesselink \paper Singularities in the  nilpotent cone of a classical group \jour Trans. Amer. Math. Soc. \vol 222 \yr 1976 \pages 1--32 \endref

\ref \no \degen  \by B. Huisgen-Zimmermann  \paper  Top-stable degenerations of finite dimensional representations I \jour Proc. London Math. Soc. \vol 96 \yr 2008 \pages 163--198 \endref

\ref \no \hierarchies \bysame \paper A hierarchy of parametrizing varieties for representations \paperinfo in Rings, Modules and Representations (N.V. Dung, et al., eds.) \jour Contemp. Math. \vol 480   \yr 2009 \pages 207--239 \endref

\ref\no \Kone \by V. Kac \paper Infinite root systems, representations of
graphs and invariant theory \jour Invent. Math. \vol 56 \yr 1980 \pages
57--92   \endref

\ref\no \Ktwo \bysame \paper Infinite root systems, representations of
graphs and invariant theory \jour J. Algebra \vol 78 \yr 1982 \pages
141--162  \endref

\ref\no \Kraft \by H.-P. Kraft \paper Geometric methods in
representation theory \inbook in Representations of Algebras (Puebla 1980) \eds M.
Auslander and E. Lluis \bookinfo Lecture Notes in Math. 944 \publ 
Spring\-er-Verlag \publaddr Berlin \yr 1982 \pages 180--258  \endref

\ref\no\Mats \by H. Matsumura \book Commutative Algebra \bookinfo Second Ed. \publaddr Reading \yr 1980 \publ Benjamin/Cummings \endref

\ref\no\MS \by E. Miller and B. Sturmfels \book Combinatorial Commutative Algebra \bookinfo Grad. Texts in Math. 227 \publaddr New York \yr 2005 \publ Springer-Verlag  \endref

\ref\no \Rie \by C. Riedtmann \paper Degenerations for representations of quivers with relations \jour Ann. Sci. Ecole Norm. Sup. \vol 19 \yr 1986 \pages 275--301 \endref

\ref\no \RRS \by C. Riedtmann, M. Rutscho, and S. O. Smal\o \paper Irreducible components of module varieties: An example \jour J. Algebra \vol 331 \yr 2011 \pages 130--144
 \endref

\ref\no \Scho \by A. Schofield \paper General representations of quivers
\jour Proc. London Math. Soc. (3) \vol 65 \yr 1992 \pages 46--64  \endref

\ref\no \Schr \by J. Schr\"oer  \paper Varieties of pairs of nilpotent matrices annihilating each other \jour Comment. Math. Helv. \vol 79 \yr 2004 \pages 396--426  \endref

\endRefs
\enddocument